\documentclass[a4paper]{amsart}
\usepackage{amsmath,amssymb,amsfonts}
\usepackage{colordvi}
\usepackage[all]{xy}

\def\dom{\backslash}
\def\mpoint{\;.}
\def\mvirg{\;,}
\def\sur{\overline}
\def\sou{\underline}

\def\mpn{\medskip\par\noindent}

\def\mmpn{\vskip 1em minus 1em\par\noindent}

\def\sp{\bigskip\par}

\def\smp{\smallskip\par}

\def\CB{{\mathcal B}}
\def\CC{{\mathcal C}}
\def\CD{{\mathcal D}}
\def\CE{{\mathcal E}}
\def\CF{{\mathcal F}}
\def\CG{{\mathcal G}}

\def\CO{{\mathcal O}}
\def\CP{{\mathcal P}}

\def\CR{{\mathcal R}}

\def\MF{\mathfrak}

\def\Ker{\operatorname{Ker}\nolimits}

\def\Im{\operatorname{Im}\nolimits}
\def\Id{\operatorname{id}\nolimits}

\def\Mod{\operatorname{Mod}\nolimits}
\def\Hom{\operatorname{Hom}\nolimits}
\def\End{\operatorname{End}\nolimits}
\def\Aut{\operatorname{Aut}\nolimits}
\def\Ext{\operatorname{Ext}\nolimits}
\def\Tor{\operatorname{Tor}\nolimits}

\def\Rad{\operatorname{Rad}\nolimits}
\def\Soc{\operatorname{Soc}\nolimits}

\def\op{^{op}}
\def\dual{^{\scriptscriptstyle\natural}}

\def\ls#1#2{{\,^{#1}\!#2}}

\def\N{\mathbb{N}}

\def\LL{\mathbb{L}}

\newcommand{\sumb}[2]{\sum_{{\scriptstyle #1}\atop {\scriptstyle #2}}}

\newcommand{\dirsum}[1]{\mathop{\bigoplus}_{#1}\limits}

\def\bigoplusl{\bigoplus\limits}

\newcommand{\edge}[2]{\xymatrix{#1\ar@{->-}[r]&#2}}
\def\marc[#1]{\ar@{-}[#1]|(.4){\object@{<}}}
\def\mard[#1]{\ar@{-}[#1]|(.5){\object@{>}}}
\def\marb[#1]{\ar@{-}[#1]|{\object+{  }}}
\newcommand{\fleche}[2]{\xymatrix@C=4ex{*!U(0.2){#1\;}&*!U(0.5){\;#2}\marc[l]}}
\newcommand{\flecheb}[2]{\xymatrix@C=4ex{*!U(0.2){#1\;}&*!U(0.1){\;#2}\marc[l]}}

\def\pf{\par\bigskip\noindent{\bf Proof~: }}
\def\endpf{~\hfill\rlap{\hspace{-1ex}\raisebox{.5ex}{\framebox[1ex]{}}\sp}\bigskip\pagebreak[3]}

\makeatletter
\renewenvironment{enumerate}{\ifnum \@enumdepth >3 \@toodeep\else
       \advance\@enumdepth \@ne
       \edef\@enumctr{enum\romannumeral\the\@enumdepth}\list
       {\csname  label\@enumctr\endcsname}{\setlength{\topsep}{1ex}
\setlength{\itemsep}{0 pt}\usecounter
         {\@enumctr}\def\makelabel##1{\hss\llap{##1}}}\fi}{\endlist}

\renewenvironment{itemize}{\ifnum \@itemdepth >3 \@toodeep\else
\advance\@itemdepth \@ne
\edef\@itemitem{labelitem\romannumeral\the\@itemdepth}
\list{\csname\@itemitem\endcsname}{\setlength{\topsep}{1ex}\setlength
{\itemsep}{0pt}\def\makelabel##1{\hss\llap{##1}}}\fi}
{\endlist}
\def\@seccntformat#1{\csname the#1\endcsname.\quad}

\def\section{\pagebreak[3]\setcounter{prop}{0}\setcounter{equation}{0}\@startsection{section}{1}{\z@}{4ex plus  6ex}{2ex}{\center\reset@font \large\bf}}
\newcommand{\subsect}[1]{\medskip\par\noindent\pagebreak[3]\refstepcounter{subsection}\refstepcounter{prop}{\bf \thesection.\arabic{prop}.\ #1.\ }}
\makeatother

\def\theprop{\thesection.\arabic{prop}}

\renewenvironment{equation}{\refstepcounter{subsection}\refstepcounter
{prop}$$}{\leqno{\bf (\theprop)}$$}

\newenvironment{rem}[1]{\refstepcounter{subsection}\refstepcounter
{prop} \mpn{{\bf \thesection.\arabic{prop}.}\ \ \bf#1.}}{\smp}
\newenvironment{enonce}[1]{\pagebreak[3]\refstepcounter{prop}\mmpn
{{\bf  \thesection.\arabic{prop}.\ #1.}}\begin{it} }{\end{it}\smp}
\def\thesection{\arabic{section}}

\newcommand{\result}[1]{\begin{enonce}{#1}}
\newcommand{\fresult}{\end{enonce}}

\begin{document}

\title[Correspondence functors and finiteness conditions]
{Correspondence functors and finiteness conditions} 

\author{Serge Bouc}
\author{Jacques Th\'evenaz}
\date\today



\begin{abstract}
We investigate the representation theory of finite sets. The correspondence functors are the functors from the category of finite sets and correspondences
to the category of $k$-modules, where $k$ is a commutative ring.
They have various specific properties which do not hold for other types of functors.
In particular, if $k$ is a field and if $F$ is a correspondence functor, then $F$ is finitely generated if and only if the dimension of $F(X)$ grows exponentially in terms of the cardinality of the finite set~$X$. Moreover, in such a case, $F$ has actually finite length. Also, if $k$ is noetherian, then any subfunctor of a finitely generated functor is finitely generated.
\end{abstract}

\maketitle


\section{Introduction}

\bigskip
\noindent
Representations of categories have been used by many authors in different contexts. The present paper is the first in a series which develops the theory in the case of the category whose objects are all finite sets and morphisms are all correspondences between finite sets.\par

For representing a category of finite sets, there are several possible choices. Pirashvili~\cite{Pi} treats the case of pointed sets and maps, while Church, Ellenberg and Farb \cite{CEF} consider the case where the morphisms are all injective maps. Putman and Sam \cite{PS} use all $k$-linear splittable injections between finite-rank free $k$-modules (where $k$ is a commutative ring). Here, we move away from such choices by using all correspondences as morphisms. The cited papers are concerned with applications to cohomological stability, while we develop our theory without any specific application in mind. The main motivation is provided by the fact that finite sets are basic objects in mathematics. Moreover, the theory turns out to have many quite surprising results, which justify the development presented here.\par

Let $\CC$ be the category of finite sets and correspondences. We define a {\em correspondence functor} over a commutative ring~$k$ to be a functor from $\CC$ to the category $k\text{-\!}\Mod$ of all $k$-modules. As much as possible, we develop the theory for an arbitrary commutative ring~$k$. However, let us start with the case when $k$ is a field. If $F$ is a correspondence functor over a field~$k$, we prove that $F$ is finitely generated if and only if the dimension of $F(X)$ grows exponentially in terms of the cardinality of the finite set~$X$ (Theorem~\ref{exponential}). In such a case, we also prove the striking fact that $F$ has finite length (Theorem~\ref{finite-length}). This result was obtained independently by Gitlin~\cite{Gi} (for a field $k$ of characteristic zero, or algebraically closed), using a criterion proved by Wiltshire-Gordon~\cite{WG}. Moreover, for finitely generated correspondence functors, we show that the Krull-Remak-Schmidt theorem holds (Proposition~\ref{KRS}) and that projective functors coincide with injective functors (Theorem~\ref{PIF}).\par

Suppose that $k$ is a field. By well-known results about representations of categories, simple correspondence functors can be classified. In our case, they are parametrized by triples $(E,R,V)$, where $E$ is a finite set, $R$ is a partial order relation on~$E$, and $V$ is a simple $k\Aut(E,R)$-module (Theorem~\ref{parametrization2}). This is the first indication of the importance of posets in our work. However, if $S_{E,R,V}$ is the simple functor parametrized by $(E,R,V)$, then it is quite hard to describe the evaluation $S_{E,R,V}(X)$ at a finite set~$X$. We will achieve this in a future paper~\cite{BT3} by giving a closed formula for its dimension.\par

A natural question when dealing with a commutative ring~$k$ is to obtain specific results when $k$ is noetherian. We follow this track in Section~\ref{Section-noetherian} and show for instance that any subfunctor of a finitely generated correspondence functor is again finitely generated (Corollary~\ref{N subseteq M}). Also, we obtain stabilization results for $\Hom$ and $\Ext$ between correspondence functors evaluated at large enough finite sets (Theorem~\ref{stabilization}).\par

This article uses essentially only standard facts from algebra and representation theory, with the following exceptions.
A few basic results in Section~\ref{Section-categories} have been imported from elsewhere, but the main exception is the algebra $\CE_E$ of essential relations on a finite set~$E$. This algebra has been analyzed in~\cite{BT1} and all its simple modules have been classified there. This uses the {\em fundamental module} $\CP_Ef_R$ associated to a finite poset $(E,R)$ and there is an explicit description of the action of relations on~$\CP_Ef_R$. All the necessary background on this algebra $\CE_E$ of essential relations is recalled in Section~\ref{Section-simple}.
It follows that our approach of the parametrization of simple functors is based on~\cite{BT1} since it uses the fundamental modules $\CP_Ef_R$ in an important way.


\section{The representation theory of categories}\label{Section-categories}

\bigskip
\noindent
Before introducing the category $\CC$ of finite sets and correspondences,
we first recall some standard facts from the representation theory of categories.
Let $\CD$ be a category and let $X$ and $Y$ be two objects of~$\CD$.
We adopt a slightly unusual notation by writing $\CD(Y,X)$ for the set of all morphisms from $X$ to~$Y$.
We reverse the order of $X$ and~$Y$ in view of having later a left action of morphisms
behaving nicely under composition.\par

We assume that $\CD$ is small (or more generally that a skeleton of~$\CD$ is small).
This allows us to talk about the {\em set} of natural transformations between two functors starting from~$\CD$.\par

Throughout this paper, $k$ denotes a commutative ring.
It will sometimes be noetherian and sometimes a field, but we shall always emphasize when we make additional assumptions.

\result{Definition} \label{linearization} The {\em $k$-linearization} of a category~$\CD$, where $k$ is any commutative ring, is defined as follows~:
\begin{itemize}
\item The objects of $k\CD$ are the objects of $\CD$.
\item For any two objects $X$ and $Y$, the set of morphisms from $X$ to~$Y$ is the free
$k$-module $k\CD(Y,X)$ with basis $\CD(Y,X)$.
\item The composition of morphisms in~$k\CD$ is the $k$-bilinear extension
$$k\CD(Z,Y) \times k\CD(Y,X) \longrightarrow k\CD(Z,X)$$
of the composition in~$\CD$.
\end{itemize}
\fresult

\result{Definition} Let $\CD$ be a category and $k$ a commutative ring.
A {\em $k$-representation} of the category~$\CD$ is a $k$-linear functor from $k\CD$
to the category $k\text{-\!}\Mod$ of $k$-modules.
\fresult

We could have defined a $k$-representation of~$\CD$ as a functor from $\CD$ to $k\text{-\!}\Mod$, but it is convenient to linearize first the category~$\CD$ (just as for group representations, where one can first introduce the group algebra).\par

If $F:k\CD\to k\text{-\!}\Mod$ is a $k$-representation of~$\CD$ and if $X$ is an object of~$\CD$,
then $F(X)$ will be called the {\em evaluation} of~$F$ at~$X$.
Morphisms in $k\CD$ act on the left on the evaluations of~$F$ by setting,
for every $m\in F(X)$ and for every morphism $\alpha\in k\CD(Y,X)$,
$$\alpha \cdot m\; :=F(\alpha)(m) \in F(Y) \mpoint$$
We often use a dot for this action of morphisms on evaluation of functors.
With our choice of notation, if $\beta\in k\CD(Z,Y)$, then
$$(\beta\alpha)\cdot m=\beta\cdot(\alpha\cdot m) \mpoint$$

The category $\CF_k(k\CD,k\text{-\!}\Mod)$ of all $k$-representations of~$\CD$ is an abelian category.
We need a small skeleton of~$\CD$ in order to have {\em sets} of natural transformations, 
which are morphisms in~$\CF_k(k\CD,k\text{-\!}\Mod)$, but we will avoid this technical discussion.
A sequence of functors 
$$0  \;\longrightarrow\; F_1  \;\longrightarrow\; F_2  \;\longrightarrow\; F_3  \;\longrightarrow\; 0$$
is exact if and only if, for every object $X$, the evaluation sequence
$$0  \;\longrightarrow\; F_1(X)  \;\longrightarrow\; F_2(X)  \;\longrightarrow\; F_3(X)  \;\longrightarrow\; 0$$
is exact.
Also, a $k$-representation of~$\CD$ is called {\em simple} if it is nonzero and has no proper nonzero subfunctor.\par

For any object $X$ of~$\CD$, consider the representable functor $k\CD(-,X)$ (which is a projective functor by Yoneda's lemma).
Its evaluation at an object $Y$ is the $k$-module $k\CD(Y,X)$,
which has a natural structure of a $(k\CD(Y,Y), k\CD(X,X))$-bimodule by composition.\par

\result{Notation}\label{LXW}
Let $X$ be an object of~$\CD$ and let $W$ be a $k\CD(X,X)$-module.
We define
$$L_{X,W}:=k\CD(-,X)\otimes_{k\CD(X,X)}W \mpoint$$
This is a $k$-representation of~$\CD$.
\fresult

This satisfies the following adjunction property.

\result{Lemma} \label{adjunction}
Let $\CF=\CF_k(k\CD,k\text{-}\Mod)$ be the category of all $k$-representations of~$\CD$
and let $X$ be an object of~$\CD$.
\begin{enumerate}
\item The functor
$$k\CD(X,X){-}\Mod \;\longrightarrow\; \CF  \;,\quad W\mapsto L_{X,W}$$
is left adjoint of the evaluation functor
$$\CF  \;\longrightarrow\; k\CD(X,X){-}\Mod \;,\quad F\mapsto F(X) \mpoint$$
In other words, for any $k$-representation $F:k\CD\to k\text{-\!}\Mod$ and any $k\CD(X,X)$-module $W$, there is a natural isomorphism
$$\Hom_{\CF}(L_{X,W}, F) \cong \Hom_{k\CD(X,X)}(W, F(X)) \mpoint$$
Moreover $L_{X,W}(X)\cong W$ as $k\CD(X,X)$-modules.
In particular, there is a $k$-algebra isomorphism $\End_\CF(L_{X,W})\cong \End_{k\CD(X,X)}(W)$.
\item The functor $k\CD(X,X){-}\Mod\longrightarrow\CF$ is right exact. It maps projective modules to projective functors,
and indecomposable modules to indecomposable functors.
\end{enumerate}
\fresult

\pf
Part (a) is straightforward and is proved in Section~2 of~\cite{Bo1}.
Part~(b) follows because this functor is left adjoint of an exact functor and satisfies the property $L_{X,W}(X)\cong W$.
\endpf

Our next result is a slight extension of the first lemma of~\cite{Bo1}.

\result{Lemma} \label{JXW}
Let $X$ be an object of~$\CD$ and let $W$ be a $k\CD(X,X)$-module.
For any object $Y$ of~$\CD$, let
$$J_{X,W}(Y):=\Big\{\sum_i\phi_i\otimes w_i \in L_{X,W}(Y) \mid
\forall\psi\in k\CD(X,Y), \sum_i (\psi\phi_i)\cdot w_i=0 \Big\} \mpoint$$
\begin{enumerate}
\item $J_{X,W}$ is the unique subfunctor of $L_{X,W}$ which is maximal
with respect to the condition that it vanishes at~$X$.
\item If $W$ is a simple module, then $J_{X,W}$ is the unique maximal subfunctor of $L_{X,W}$
and $L_{X,W}/J_{X,W}$ is a simple functor.
\end{enumerate}
\fresult

\pf
The proof is sketched in Lemma~2.3 of~\cite{BST} in the special case of biset functors for finite groups,
but it extends without change to representations of an arbitrary category~$\CD$.
\endpf

Lemma~\ref{JXW} is our main tool for dealing with simple functors. We first fix the notation.

\result{Notation}\label{SXW}
Let $X$ be an object of~$\CD$ and let $W$ be a $k\CD(X,X)$-module.
We define
$$S_{X,W}:=L_{X,W}/J_{X,W} \mpoint$$
If $W$ is a simple $k\CD(X,X)$-module, then $S_{X,W}$ is a simple functor.
\fresult

We emphasize that $L_{X,W}$ and $S_{X,W}$ are defined for any $k\CD(X,X)$-module~$W$ and any commutative ring~$k$.
Note that we always have $J_{X,W}(X)=\{0\}$ because if $a=\sum_i\limits\phi_i\otimes w_i \in J_{X,W}(X)$,
then $a=\Id_X\otimes(\sum_i\limits\phi_i\cdot w_i) =0$.

Therefore, we have isomorphisms of $k\CD(X,X)$-modules
$$L_{X,W}(X)\cong S_{X,W}(X)\cong W \mpoint$$

\result{Proposition} \label{generated}
Let $S$ be a simple $k$-representation of~$\CD$ and let $Y$ be an object of~$\CD$ such that $S(Y)\neq 0$.
\begin{enumerate}
\item $S(Y)$ is a simple $k\CD(Y,Y)$-module.
\item $S\cong S_{Y,S(Y)}$.
\item $S$ is generated by $S(Y)$, that is, $S(X)=k\CD(X,Y)S(Y)$ for all objects~$X$.
More precisely, if $0\neq u\in S(Y)$, then $S(X)=k\CD(X,Y)\cdot u$.
\end{enumerate}
\fresult

\pf
(c) Given $0\neq u\in S(Y)$, let $S'(X)=k\CD(X,Y)\cdot u$ for all objects~$X$.
This clearly defines a nonzero subfunctor $S'$ of~$S$, so $S'=S$ by simplicity of~$S$.\mpn

(a) This follows from (c).\mpn

(b) By the adjunction of Lemma~\ref{adjunction}, the identity $\Id:S(Y)\to S(Y)$ corresponds to a non-zero morphism $\theta:L_{Y,S(Y)}\to S$. Since $S$ is simple, $\theta$ must be surjective. But $S_{Y,S(Y)}$ is the unique simple quotient of $L_{Y,S(Y)}$, by Lemma~\ref{JXW} and Notation~\ref{SXW}, so $S\cong S_{Y,S(Y)}$.
\endpf

It should be noted that a simple $k$-representation $S$ has many possible realizations $S\cong S_{Y,W}$ as above, where $W=S(Y)\neq0$.
However, if there is a notion of unique minimal object, then one can parametrize simple functors $S$ by setting $S\cong S_{Y,W}$, where $Y$ is the unique minimal object such that $S(Y)\neq0$ (see Theorem~\ref{parametrization1} for the case of correspondence functors).\par

Our next proposition is Proposition~3.5 in~\cite{BST} in the case of biset functors,
but it holds more generally and we just recall the proof of~\cite{BST}.

\result{Proposition} \label{subquotients}
Let $S$ be a simple $k$-representation of~$\CD$ and let $Y$ be an object of~$\CD$ such that $S(Y)\neq 0$.
Let $F$ be any $k$-representation of~$\CD$. Then the following are equivalent:
\begin{enumerate}
\item $S$ is isomorphic to a subquotient of~$F$.
\item The simple $k\CD(Y,Y)$-module $S(Y)$ is isomorphic to a subquotient of the $k\CD(Y,Y)$-module $F(Y)$.
\end{enumerate}
\fresult

\pf
It is clear that (a) implies (b). Suppose that (b) holds and let $W_1$, $W_2$ be submodules of~$F(Y)$ such that $W_2\subset W_1$ and $W_1/W_2\cong S(Y)$. For $i\in\{1,2\}$, let $F_i$ be the subfunctor of $F$ generated by $W_i$. Explicitly, for any object $X$ of~$\CD$, 
$F_i(X)=k\CD(X,Y)\cdot W_i\subseteq F(X)$. Then $F_i(Y)=W_i$ and $(F_1/F_2)(Y)=W_1/W_2\cong S(Y)$. The isomorphism $S(Y)\to (F_1/F_2)(Y)$ induces, by the adjunction of Lemma~\ref{adjunction}, a nonzero morphism $\theta: L_{Y,S(Y)}\to F_1/F_2$. Since $S(Y)$ is simple, $L_{Y,S(Y)}$ has a unique maximal subfunctor $J_{Y,S(Y)}$, by Lemma~\ref{JXW}, and $L_{Y,S(Y)}/J_{Y,S(Y)} \cong S_{Y,S(Y)} \cong S$, by Proposition~\ref{generated}. Let $F'_1=\theta(L_{Y,S(Y)})$ and $F'_2=\theta(J_{Y,S(Y)})$. Since $\theta\neq0$, we obtain
$$F'_1/F'_2 \cong L_{Y,S(Y)}/J_{Y,S(Y)} \cong S_{Y,S(Y)} \cong S \,,$$
showing that $S$ is isomorphic to a subquotient of~$F$.\par

(Actually, as observed by Hida and Yagita in Lemma~3.1 of~\cite{HY}, we have an equality $F'_1=F_1$, because both subfunctors are generated by their commun evaluation at~$Y$.)
\endpf


\section{Correspondence functors}\label{Section-functors}

\bigskip
\noindent
Leaving the general case, we now prepare the ground for the category $\CC$ we are going to work with.

\result{Definition} Let $X$ and $Y$ be sets.
\begin{enumerate}
\item A {\em correspondence} from $X$ to~$Y$ is a subset of the cartesian product~$Y\times X$.
Note that we have reversed the order of $X$ and~$Y$ for the reasons mentioned at the beginning of Section~\ref{Section-categories}.
\item A correspondence is often called a relation
but we use this terminology only when $X=Y$, in which case we say that a subset of $X\times X$ is a {\em relation on~$X$}.
\item If $\sigma$ is a permutation of~$X$, then there is a corresponding relation on~$X$ which we write
$$\Delta_\sigma:=\{(\sigma(x),x)\in X\times X\mid x\in X\}  \mpoint$$
In particular, when $\sigma=\Id$, we also write
$$\Delta_{\Id}=\Delta_X=\{(x,x)\in X\times X\mid x\in X\}\mpoint$$
\end{enumerate}
\fresult

\result{Definition} Let $\CC$ denote the following category~:
\begin{itemize}
\item The objects of $\CC$ are the finite sets.
\item For any two finite sets $X$ and $Y$, the set $\CC(Y,X)$ is the set of all correspondences from $X$ to~$Y$.
\item The composition of correspondences is as follows. Given $R\subseteq Z\times Y$ and $S\subseteq Y\times X$,
then $RS$ is defined by
$$RS:=\{ \, (z,x)\in Z\times X \,\mid\, \exists\;y\in Y \;\text{ such that } \; (z,y)\in R \,\text{ and }\, (y,x)\in S \,\} \mpoint$$
\end{itemize}
\fresult

The identity morphism $\Id_X$ is the diagonal subset $\Delta_X\subseteq X\times X$
(in other words the equality relation on~$X$).

\result{Definition} Let $k\CC$ be the linearization of the category $\CC$, where $k$ is any commutative ring (see Definition~\ref{linearization}).
\begin{enumerate}
\item A {\em correspondence functor (over~$k$)} is a $k$-representation of the category~$\CC$,
that is, a $k$-linear functor from $k\CC$ to the category $k\hbox{\rm-Mod}$ of $k$-modules.
\item We let $\CF_k=\CF_k(k\CC,k\hbox{\rm -Mod})$ be the category of all such correspondence functors (an abelian category).
\end{enumerate}
\fresult

In part (b), we need to restrict to a small skeleton of~$\CC$ in order to have {\em sets} of natural transformations, which are morphisms in~$\CF_k$,
but we avoid this technical discussion.
It is clear that $\CC$ has a small skeleton, for instance by taking the full subcategory having one object for each cardinality.\par

For any finite set $E$, we define
$$\CR_E:=k\CC(E,E) \mvirg$$
the $k$-algebra of the monoid $\CC(E,E)$ of all relations on~$E$, in other words the algebra of the semigroup of Boolean matrices of size~$|E|$.
The representable functor $k\CC(-,E)$ (sometimes called Yoneda functor)
is the very first example of a correspondence functor.
By definition, it is actually isomorphic to the functor $L_{E,\CR_E}$, see Notation~\ref{LXW}.
If $W$ is an $\CR_E$-module generated by a single element~$w$ (for instance a simple module),
then the functor $L_{E,W}$ is isomorphic to a quotient of $k\CC(-,E)$ via the surjective homomorphism
$$k\CC(-,E) \;\longrightarrow\; L_{E,W} = k\CC(-,E) \otimes_{\CR_E} W \,,
\qquad \phi \mapsto \phi\otimes w \mpoint$$
The representable functor $k\CC(-,E)$ is projective (by Yoneda's lemma).\par

Our next result is basic, but has several important corollaries.

\result{Lemma} \label{summand1} Let $E$ and $F$ be finite sets with $|E| \leq |F|$.
There exist correspondences $i_*\in \CC(F,E)$ and $i^*\in \CC(E,F)$ such that $i^*i_*=\Id_E$.
\fresult

\pf
Since $|E| \leq |F|$, there exists an injective map $i:E\hookrightarrow F$.
Let $i_*\subseteq F\times E$ denote the correspondence 
$$i_*=\big\{\big(i(e),e\big)\mid e\in E\big\}\mvirg$$
and $i^*\subseteq E\times F$ denote the correspondence 
$$i^*=\big\{\big(e,i(e)\big)\mid e\in E\big\}\mpoint$$
As $i$ is injective, one checks easily that $i^*i_*=\Delta_E$, that is, $i^*i_*=\Id_E$.
\endpf

In other words, this lemma says that the object $E$ of~$\CC$ behaves like a direct summand of the object~$F$ whenever $|E|\leq|F|$.

\result{Corollary} \label{summand2} Let $E$ and $F$ be finite sets with $|E| \leq |F|$.
The representable functor $k\CC(-,E)$ is isomorphic to a direct summand of the representable functor~$k\CC(-,F)$.
\fresult

\pf
Right multiplication by~$i^*$ defines a homomorphism of correspondence functors
$$k\CC(-,E) \longrightarrow k\CC(-,F) \mvirg$$
and right multiplication by~$i_*$ defines a homomorphism of correspondence functors
$$k\CC(-,F) \longrightarrow k\CC(-,E) \mpoint$$
Their composite is the identity of $k\CC(-,E)$, because $i^*i_*=\Id_E$.
\endpf

\result{Corollary} \label{projective-CEF} Let $E$ and $F$ be finite sets with $|E| \leq |F|$.
The left $k\CC(F,F)$-module $k\CC(F,E)$ is projective.
\fresult

\pf
By Corollary~\ref{summand2}, $k\CC(F,E)$ is isomorphic to a direct summand of~$k\CC(F,F)$, which is free.
\endpf

\result{Corollary} \label{vanishing-below} Let $E$ and $F$ be finite sets with $|E| \leq |F|$.
Let $M$ be a correspondence functor. If $M(F)=0$, then $M(E)=0$.
\fresult

\pf
For any $m\in M(E)$, we have $m=i^*i_*\cdot m$. But $i_*\cdot m\in M(F)$, so $i_*\cdot m=0$.
Therefore $m=0$.
\endpf

\result{Corollary} \label{composition-iso} Let $E$ and $F$ be finite sets with $|E| \leq |F|$.
For every finite set~$X$, composition in the category $k\CC$
$$\mu: k\CC(X,F)\otimes_{\CR_F}k\CC(F,E) \longrightarrow k\CC(X,E)$$
is an isomorphism.
\fresult

\pf
The inverse of $\mu$ is given by
$$\phi:k\CC(X,E)\longrightarrow k\CC(X,F)\otimes_{\CR_F}k\CC(F,E) \,,\qquad
\phi(\alpha)= \alpha i^* \otimes i_* \mpoint$$
Composing with $\mu$, we obtain $\mu\phi(\alpha)=\mu(\alpha i^* \otimes i_*)= \alpha i^* i_*=\alpha$, so $\mu\phi=\Id$.
On the other hand, if $\beta\in k\CC(X,F)$ and $\gamma\in k\CC(F,E)$, then
$\gamma i^*$ belongs to $\CR_F$ and therefore
$$\phi\mu(\beta\otimes\gamma)=\phi(\beta\gamma)=\beta\gamma i^* \otimes i_*= \beta\otimes\gamma i^* i_* =\beta\otimes\gamma \mvirg$$
showing that $\phi\mu=\Id$.
\endpf

Now we move to direct summands of representable functors, given by some idempotent.
If $R$ is an idempotent in~$\CR_E$, then $k\CC(-,E)R$ is a direct summand of~$k\CC(-,E)$,
hence projective again. In particular, if $R$ is a preorder on~$E$, that is, a relation which is reflexive and transitive,
then $R$ is idempotent (because $R\subseteq R^2$ by reflexivity and $R^2\subseteq R$ by transitivity).
There is an equivalence relation~$\sim$ associated with~$R$, defined by
$$x\sim y \;\Longleftrightarrow\;  (x,y)\in R \;\text{ and }\; (y,x)\in R \mpoint$$
Then $R$ induces an order relation~$\overline R$ on the quotient set $\overline E=E/\sim$ such that
$$(x,y)\in R \;\Longleftrightarrow\;  (\overline x,\overline y)\in \overline R \mvirg$$
where $\overline x$ denotes the equivalence class of~$x$ under~$\sim$.

\result{Lemma} \label{preorder}
Let $E$ be a finite set and let $R$ be a preorder on~$E$. Let $\overline R$ be the corresponding order
on the quotient set $\overline E=E/\sim$ and write $e\mapsto \overline e$ for the quotient map $E\to \overline E$.
The correspondence functors $k\CC(-,E)R$ and $k\CC(-,\overline E)\overline R$
are isomorphic via the isomorphism
$$k\CC(-,E)R \;\longrightarrow\; k\CC(-,\overline E)\overline R\,,
\qquad S \mapsto \overline S \mvirg$$
where for any correspondence $S\subseteq X\times E$,
the correspondence $\overline S\subseteq X\times \overline E$ is defined by
$$(x,\overline e)\in \overline S \;\Longleftrightarrow\;  (x,e) \in S \mpoint$$
\fresult

\pf It is straightforward to check that $\overline S$ is well-defined.
If $S \in\CC(X,E)R$, then $S=SR$ and it follows that $\overline S=\overline S\,\overline R$.
Surjectivity is easy and injectivity follows from the definition of~$\overline S$.
\endpf

This shows that it is relevant to consider the functors $k\CC(-,E)R$ where $R$ is an order on~$E$.
We will see later that they play an important role in connection with simple functors.\par

We end this section with the definition of duality.

\result{Proposition} Let $X$ and $Y$ be finite sets. If $R\subseteq Y\times X$, let $R\op$ denote the {\em opposite correspondence}, defined by
$$R\op:=\{(x,y)\in X\times Y\mid (y,x)\in R\}\subseteq X\times Y\mpoint$$
Then the assignment $R\mapsto R\op$ induces an isomorphism from $\CC$ to the opposite category $\CC\op$, which extends to an isomorphism from $k\CC$ to $k\CC\op$.
\fresult

\pf Let $X$, $Y$, $Z$ be finite sets, $R\subseteq Y\times X$ and $S\subseteq Z\times Y$. One checks easily that $(SR)\op=R\op S\op$.\endpf

We use opposite correspondences to define dual functors. The notion will be used in Section~\ref{Section-projective}.
 
\result{Definition}\label{def-dual} Let $F$ be a correspondence functor over $k$.
The {\em dual} $F\dual$ of $F$ is the correspondence functor defined on a finite set $X$ by
$$F\dual(X):=\Hom_k\big(F(X),k\big)\mpoint$$
If $Y$ is a finite set and $R\subseteq Y\times X$, then the map $F\dual(R): F\dual(X)\to F\dual(Y)$ is defined by
$$\forall \alpha\in F\dual(X),\;\;F\dual(R)(\alpha):=\alpha\circ F(R\op)\mpoint$$
\fresult


\section{The parametrization of simple correspondence functors}\label{Section-simple}

\bigskip
\noindent
In order to study simple modules or simple functors, it suffices to work over a field~$k$,
by standard commutative algebra. If we assume that $k$ is a field, then
the evaluation at a finite set $X$ of a representable functor, or of a simple functor,
is always a finite-dimensional $k$-vector space.
As before, we continue to work with an arbitrary commutative base ring~$k$ and assume that it is a field when necessary.

\result{Definition} A {\em minimal set} for a correspondence functor~$F$ is a finite set~$X$ of minimal cardinality such that $F(X)\neq 0$.
For a nonzero functor, such a minimal set always exists and is unique up to bijection.
\fresult

Our next task is to describe the parametrization of simple correspondence functors.
This uses the algebra
$$\CR_X:=k\CC(X,X)$$
of all relations on~$X$.
This algebra was studied in~\cite{BT1} and we use this approach.
A relation~$R$ on~$X$ is called {\em essential\/} if it does not factor through a set of cardinality strictly smaller than~$|X|$. In other words, $R$ has maximal Schein rank in the sense of~Section~1.4 of~\cite{Ki}.
The $k$-submodule generated by the set of inessential relations is a two-sided ideal
$$I_X:=\sum_{|Y|<|X|} k\CC(X,Y)k\CC(Y,X)$$
and the quotient
$$\CE_X:=k\CC(X,X)/I_X$$
is called the {\em essential algebra}.
A large part of its structure has been elucidated in~\cite{BT1}.\par

The following parametrization theorem is similar to Proposition~2 in~\cite{Bo1} or Theorem~4.3.10 in~\cite{Bo2}.
The context here is different, but the proof is essentially the same.

\result{Theorem} \label{parametrization1} Assume that $k$ is a field.
\begin{enumerate}
\item Let $S$ be a simple correspondence functor,
let $E$ be a minimal set for~$S$, and let $W=S(E)$.
Then $W$ is a simple module for the essential algebra $\CE_E$ (with $I_E$ acting by zero) and $S\cong S_{E,W}$.
\item Let $E$ be a finite set and let $W$ be a simple module for the essential algebra $\CE_E$,
viewed as a module for the algebra $\CR_E$ by making $I_E$ act by zero on~$W$.
Then $E$ is a minimal set for $S_{E,W}$.
Moreover, $S_{E,W}(E)\cong W$ (as $\CE_E$-modules).
\item The set of isomorphism classes of simple correspondence functors is par\-ametrized by
the set of isomorphism classes of pairs $(E,W)$
where $E$ is a finite set and $W$ is a simple $\CE_E$-module.
\end{enumerate}
\fresult

\pf
(a) Since $S(Y)=\{0\}$ if $|Y|<|E|$, we have
$$I_E\cdot S(E)=\sum_{|Y|<|E|} k\CC(E,Y)k\CC(Y,E)\cdot S(E)\subseteq \sum_{|Y|<|E|} k\CC(E,Y) \cdot S(Y) =\{0\} \mvirg$$
so $S(E)$ is a module for the essential algebra $\CE_E$.
Now the identity of~$S(E)$ corresponds by adjunction
to a nonzero homomorphism $L_{E,W}\to S$, where $W=S(E)$ (see Lemma~\ref{adjunction}).
This homomorphism is surjective since $S$ is simple.
But $L_{E,W}$ has a unique simple quotient, namely $S_{E,W}$, hence $S\cong S_{E,W}$.\mpn

(b) Suppose that $S_{E,W}(Y)\neq\{0\}$. Then $L_{E,W}(Y)\neq J_{E,W}(Y)$,
so there exists a correspondence $\phi\in \CC(Y,E)$ and $v\in W$ such that $\phi\otimes v\in L_{E,W}(Y)-J_{E,W}(Y)$.
By definition of~$J_{E,W}$, this means that there exists a correspondence $\psi\in \CC(E,Y)$
such that $\psi\phi\cdot v\neq 0$. Since $W$ is a module for the essential algebra $\CE_E=\CR_E/I_E$,
we have $\psi\phi\notin I_E$. But $\psi\phi$ factorizes through~$Y$, so we must have $|Y|\geq |E|$.
Thus $E$ is a minimal set for $S_{E,W}$. The isomorphism $S_{E,W}(E)\cong W$ is a general fact mentioned before.\mpn

(c) This follows from (a) and~(b).
\endpf

Theorem~\ref{parametrization1} reduces the classification of simple correspondence functors
to the question of classifying all simple modules for the essential algebra $\CE_E$. 
Fortunately, this has been achieved in~\cite{BT1}.
An alternative path would be to use the classical approach to simple $\CR_E$-modules via idempotents in the semigroup $\CC(E,E)$ and Green's theory of $J$-classes (see the textbook~\cite{CP}, or the recent article \cite{GMS} for a modern point of view).
We do not follow this approach because we need later in an important way a fundamental module, defined below, which is not part of the classical approach, but is the key for the classification in~\cite{BT1}.\par

The simple $\CE_E$-modules are actually modules for a quotient $\CP_E=\CE_E/N$ where $N$ is a nilpotent ideal defined in~\cite{BT1}. We call $\CP_E$ the {\em algebra of permuted orders},
because it has a $k$-basis consisting of all relations on~$E$ of the form
$\Delta_\sigma R$, where $\sigma$ runs through the symmetric group $\Sigma_E$ of all permutations of~$E$,
and $R$ is an order on~$E$. By an order, we always mean a partial order relation.
The product of two orders $R$ and $S$ in~$\CP_E$ is the transitive closure of $R\cup S$ if this transitive closure is an order, and zero otherwise.\par

We let $\CO$ be the set of all orders on~$E$
and $\Aut(E,R)$ the stabilizer of the order~$R$ under the action of the symmetric group~$\Sigma_E$ by conjugation.
For the description of simple $\CE_E$-modules,
we need the following new basis of~$\CP_E$ (see Theorem~6.2 in~\cite{BT1} for details).

\result{Lemma} \label{idempotents}
\begin{enumerate}
\item There is a set $\{f_R\mid R\in\CO\}$ of orthogonal idempotents of~$\CP_E$ whose sum is~1, such that
$\CP_E$ has a $k$-basis consisting of all elements of the form
$\Delta_\sigma f_R$, where $\sigma\in \Sigma_E$ and $R\in\CO$.
\item For any $\sigma\in\Sigma_E$, we have $\ls\sigma f_R=f_{\ls\sigma R}$, where $\ls\sigma x=\Delta_\sigma x\Delta_{\sigma^{-1}}$ for any $x\in \CP_E$. In particular, $\Delta_\sigma f_R\Delta_{\sigma^{-1}}=f_R$ if $\sigma\in\Aut(E,R)$.
\item For any order~$Q$ on~$E$, we have~:
$$Qf_R \neq 0 \iff Qf_R=f_R \iff Q\subseteq R \mpoint$$
\end{enumerate}
\fresult

\medskip
For the description of simple $\CE_E$-modules and then simple correspondence functors,
we will make use of the left $\CE_E$-module $\CP_E f_R$.
This module is actually defined without assuming that $k$ is a field.

\result{Definition}
Let $(E,R)$ be a finite poset (i.e. $E$ is a finite set and $R$ is an order on~$E$).
We call $\CP_E f_R$ the {\em fundamental module} for the algebra~$\CE_E$, associated with the poset $(E,R)$.
\fresult

We now describe its structure.

\result{Proposition} \label{fundamental-module}
Let $(E,R)$ be a finite poset
\begin{enumerate}
\item The fundamental module $\CP_E f_R$ is a left module for the algebra $\CP_E$, hence also a left module for the essential algebra $\CE_E$ and for the algebra of relations~$\CR_E$.
\item $\CP_E f_R$ is a free $k$-module with a $k$-basis consisting of the elements $\Delta_\sigma f_R$,
where $\sigma$ runs through the group $\Sigma_E$ of all permutations of $E$.
\item $\CP_E f_R$ is a $(\CP_E,k\Aut(E,R))$-bimodule and the right action of $k\Aut(E,R)$ is free.
\item The action of the algebra of relations $\CR_E$ on the module $\CP_E f_R$
is given as follows. For any relation $Q\in \CC(E,E)$,
$$Q\cdot\Delta_\sigma f_R=\left\{\begin{array}{ll}
\Delta_{\tau\sigma}f_R&\hbox{if}\;\;\exists\tau\in\Sigma_E\;\hbox{such that}\;
\Delta_E\subseteq \Delta_{\tau^{-1}}Q\subseteq {\ls\sigma R},\\
0&\hbox{otherwise}\mvirg\end{array}\right.$$
where $\Delta_E$ is the diagonal of $E\times E$, and $\ls\sigma R=\big\{\big(\sigma(e),\sigma(f)\big)\mid (e,f)\in R\big\}$ (recall that $\tau$ is unique in the first case).
\end{enumerate}
\fresult

\pf
See Corollary 7.3 and Proposition~8.5 in~\cite{BT1}.
\endpf

The description of all simple $\CE_E$-modules is as follows (see Theorem~8.1 in~\cite{BT1} for details).

\result{Theorem} \label{simple-modules} Assume that $k$ is a field and let $E$ be a finite set.
\begin{enumerate}
\item Let $R$ be an order on~$E$ and let $\CP_E f_R$ be the corresponding fundamental module.
If $V$ is a simple $k\Aut(E,R)$-module, then
$$T_{R,V}:=\CP_E f_R\otimes_{k\Aut(E,R)} V$$
is a simple $\CP_E$-module (hence also a simple $\CE_E$-module).
\item Every simple $\CE_E$-module is isomorphic to a module $T_{R,V}$ as in~(a).
\item For any permutation $\sigma\in \Sigma_E$, we have $T_{\ls\sigma R,\ls\sigma V}\cong T_{R,V}$,
where $\ls\sigma R=\Delta_\sigma R\Delta_{\sigma^{-1}}$ is the conjugate order
and $\ls\sigma V$ is the conjugate module.
\item The set of isomorphism classes of simple $\CE_E$-modules is parametrized by
the set of conjugacy classes of pairs $(R,V)$ where $R$ is an order on~$E$ and $V$ is a simple $k\Aut(E,R)$-module.
\end{enumerate}
\fresult

Putting together Theorem~\ref{parametrization1} and Theorem~\ref{simple-modules}, we finally obtain the following parametrization, which is essential for our purposes.

\result{Theorem} \label{parametrization2} Assume that $k$ is a field.
The set of isomorphism classes of simple correspondence functors is parametrized by
the set of isomorphism classes of triples $(E,R,V)$
where $E$ is a finite set, $R$ is an order on~$E$, and $V$ is a simple $k\Aut(E,R)$-module.
\fresult

\pf By Theorem~\ref{parametrization1}, the isomorphism classes of simple correspondence functors are parametrized by pairs $(E,W)$,
where $E$ is a finite set and $W$ is a simple $\CE_E$-module.
Now by Theorem~\ref{parametrization2}, the isomorphism classes of simple $\CE_E$-modules~$W$ are parametrized by pairs $(R,V)$,
where $R$ is an order on~$E$ and $V$ is a simple $k\Aut(E,R)$-module.
Hence the pairs $(E,W)$ become triples $(E,R,V)$.
\endpf

\result{Notation} \label{notation-SERV}   
Let $(E,R)$ be a finite poset.
\begin{enumerate}
\item If $V$ is a simple $k\Aut(E,R)$-module,
we denote by $S_{E,R,V}$ the simple correspondence functor parametrized by the triple $(E,R,V)$.
Explicitly,
$$S_{E,R,V}=L_{E,T_{R,V}}/J_{E,T_{R,V}} \mpoint$$
\item More generally, for any commutative ring $k$ and any $k\Aut(E,R)$-module~$V$, we define
the $\CP_E$-module (hence also an $\CR_E$-module)
$$T_{R,V}:=\CP_E f_R\otimes_{k\Aut(E,R)} V$$
and the correspondence functor
$$S_{E,R,V}:=L_{E,T_{R,V}}/J_{E,T_{R,V}} \mpoint$$
\end{enumerate}
\fresult

We end this section with a basic result concerning the correspondence functors $S_{E,R,V}$,
where $k$ is any commutative ring and $V$ is any $k\Aut(E,R)$-module.

\result{Lemma}\label{SERV}
Let $(E,R)$ be a finite poset and let $V$ be any $k\Aut(E,R)$-module.
\begin{enumerate}
\item $E$ is a minimal set for~$S_{E,R,V}$.
\item $S_{E,R,V}(E)\cong T_{R,V}$ as $\CR_E$-modules.
\end{enumerate}
\fresult

\pf
Both $\CP_E f_R$ and $T_{R,V}=\CP_E f_R\otimes_{k\Aut(E,R)} V$ are left modules for the essential algebra $\CE_E$.
Therefore, the argument given in part~(b) of Theorem~\ref{parametrization1} shows again that
$E$ is a minimal set for~$S_{E,R,V}$.
Moreover, since $J_{E,T_{R,V}}$ vanishes on evaluation at~$E$, we have
$$S_{E,R,V}(E)\cong L_{E,T_{R,V}}(E)=T_{R,V}=\CP_E f_R\otimes_{k\Aut(E,R)} V \mvirg$$
as required.
\endpf

When $k$ is a field and $V$ is simple, we recover the simple module $S_{E,R,V}(E)=T_{R,V}$ of Theorem~\ref{simple-modules}.


\section{Small examples} \label{Section-examples}

\bigskip
\noindent
In this section, we describe two small examples.

\subsect{Example} \label{empty}
Let $E=\emptyset$ be the empty set and consider the representable functor $k\CC(-,\emptyset)$.
Then $\CC(X,\emptyset)=\{\emptyset\}$ is a singleton for any finite set~$X$,
so $k\CC(X,\emptyset)\cong k$. Moreover, any correspondence $S\in\CC(Y,X)$ is sent to the identity map
from $k\CC(X,\emptyset)\cong k$ to $k\CC(Y,\emptyset)\cong k$.
This functor deserves to be called the {\em constant functor}. We will denote it by $\sou{k}$.\par

Assume that $k$ is a field. The algebra $k\CC(\emptyset,\emptyset)\cong k$ has a unique simple module~$k$.
It is then easy to check that $L_{\emptyset,k}=k\CC(-,\emptyset)$ and $J_{\emptyset,k}=\{0\}$.
Therefore$$k\CC(-,\emptyset) = S_{\emptyset, \emptyset,k}$$
the simple functor indexed by $(\emptyset, \emptyset,k)$.
Here the second $\emptyset$ denotes the only relation on the empty set,
while $k$ is the only simple module for the group algebra
$k\Aut(E,R)=k\Sigma_\emptyset\cong k$, where $\Sigma_\emptyset=\{\Id\}$ is the symmetric group of the empty set.
Note that $S_{\emptyset, \emptyset,k}$ is projective, because it is a representable functor. 

\subsect{Example} \label{one}
Let $E=\{1\}$ be a set with one element and consider the representable functor $\CC(-,\{1\})$.
Then $\CC(X,\{1\})$ is in bijection with the set $\CB(X)$ of all subsets of~$X$,
because $X\times \{1\} \cong X$. It is easy to see that a correspondence $S\in\CC(Y,X)$ is sent to the map
$$\CB(X) \longrightarrow \CB(Y) \,,\qquad A \mapsto S_A \mvirg$$
where $S_A=\{ y\in Y \mid \exists\, a\in A \,\text{ such that } \, (y,a)\in S \}$.
Thus $k\CB\cong k\CC(-,\{1\})$ is a correspondence functor such that $k\CB(X)$ is a free $k$-module
with basis~$\CB(X)$ and rank $2^{|X|}$ for every finite set~$X$.\par

The functor $k\CB$ has a subfunctor isomorphic to the constant functor $S_{\emptyset, \emptyset,k}$,
because $\CB(X)$ contains the element $\emptyset$ which is mapped to $\emptyset$ by any correspondence.
We claim that, if $k$ is a field, the quotient $k\CB/S_{\emptyset, \emptyset,k}$ is a simple functor.\par

Assume that $k$ is a field. The algebra $k\CC(\{1\},\{1\})$ has dimension~2,
actually isomorphic to $k\times k$ with two primitive idempotents $\emptyset$ and $\{(1,1)\}-\emptyset$.
The essential algebra $\CE_{\{1\}}$ is a one-dimensional quotient
and its unique simple module~$W$ is one-dimensional and corresponds to the pair $(R,k)$,
where $R$ is the only order relation on~$\{1\}$ and $k$ is the only simple module for the group algebra
$k\Aut(E,R)=k\Sigma_{\{1\}}\cong k$, with $\Sigma_{\{1\}}=\{\Id\}$ the symmetric group of~$\{1\}$.
Thus there is a simple functor $S_{\{1\},W} = S_{\{1\},R,k}$.\par

The kernel of the quotient map
$$k\CB\cong k\CC(-,\{1\}) \;\longrightarrow \; k\CC(-,\{1\})\otimes_{k\CC(\{1\},\{1\})} W  =L_{\{1\},W}$$
is the constant subfunctor $S_{\emptyset, \emptyset,k}$ mentioned above,
because $\emptyset \in \CC(X,\{1\})$ can be written
$\emptyset \cdot\emptyset$, with the second empty set belonging to $\CC(\{1\},\{1\})$, thus acting by zero on~$W$.
Now we know that $L_{\{1\},W}/J_{\{1\},W} = S_{\{1\},W}$ and we are going to show that $J_{\{1\},W}=\{0\}$.
It then follows that $L_{\{1\},W} = S_{\{1\},W}$ is simple and isomorphic to $k\CB/S_{\emptyset, \emptyset,k}$,
proving the claim above.\par

In order to prove that $J_{\{1\},W}=\{0\}$, we let $u\in L_{\{1\},W}(X)$, which can be written
$$u=\sum_{A\in \CB(X)} \lambda_A (A\times\{1\}) \otimes w\mvirg$$
where $w$ is a generator of~$W$ and $\lambda_A\in k$ for all $A$.
Since the empty set acts by zero on~$w$, the sum actually runs over nonempty subsets $A\in\CB(A)$.
Then $u\in J_{\{1\},W}(X)$ if and only if, for all $(\{1\}\times B)\in\CC(\{1\},X)$,
we have
$$\sum_{A\in \CB(X)} \lambda_A \; (\{1\}\times B)(A\times\{1\}) \cdot w=0 \mpoint$$
Since $\emptyset$ acts by zero and $\{1\}$ acts as the identity, we obtain 
$$(\{1\}\times B)(A\times\{1\}) \cdot w= \left\{ 
\begin{array}{rl} 0 \quad &\text{if } \; B\cap A=\emptyset \,, \\
w \quad &\text{if } \; B\cap A\neq \emptyset \mpoint
\end{array} \right. $$
This yields the condition
$$\sum_{A\cap B\neq\emptyset} \lambda_A=0 \,, \qquad \text{for every nonempty } \;B\in\CB(X) \mpoint$$
We prove by induction that $\lambda_C=0$ for every nonempty $C\in \CB(X)$.
Subtracting the condition for $B=X$ and for $B=X-C$, we obtain
$$0=\sum_{A\neq\emptyset} \lambda_A - \sum_{A\not\subseteq C} \lambda_A
=\sum_{\emptyset\neq A \subseteq C} \lambda_A \mpoint$$
If $C=\{c\}$ is a singleton, we obtain $\lambda_{\{c\}}=0$ and this starts the induction.
In the general case, we obtain by induction $\lambda_A=0$ for $\emptyset\neq A \neq C$,
so we are left with $\lambda_C=0$. Therefore $u=0$ and so $J_{\{1\},W}=\{0\}$.\par

There is a special feature of this small example, namely that the exact sequence
$$0 \longrightarrow S_{\emptyset, \emptyset,k} \longrightarrow k\CB \longrightarrow S_{\{1\},R,k}
\longrightarrow 0 $$
splits. This is because there is a retraction $k\CB \to S_{\emptyset, \emptyset,k}$ defined by
$$k\CB(X) \longrightarrow k \,,\qquad A\mapsto 1 \mvirg$$
which is easily checked to be a homomorphism of functors. Since $k\CB$ is projective
(because it is a representable functor), its direct summand $S_{\{1\},R,k}$ is projective.

\begin{rem}{Remark}
In both Example \ref{empty} and Example \ref{one},
there is a unique order relation~$R$ on~$E$, which is a total order.
Actually, these examples are special cases of the general situation of a total order,
which is studied in~\cite{BT2}.
\end{rem}


\section{Finite generation} \label{Section-fg}

\bigskip
\noindent
In this section, we analyze the property of finite generation for correspondence functors.

\result{Definition} Let $\{E_i\mid i\in I\}$ be a family of finite sets indexed by a set~$I$
and, for every $i\in I$, let $m_i\in M(E_i)$.
A correspondence functor $M$ is said to be {\em generated} by the set $\{m_i\mid i\in I\}$ if for every finite set~$X$
and every element $m\in M(X)$, there exists a finite subset $J\subseteq I$ such that
$$m=\sum_{j\in J} \alpha_jm_j \,, \qquad \text{for some } \alpha_j\in k\CC(X,E_j) \quad(\text{where }j\in J) \mpoint$$
In the case where $I$ is finite, $M$ is said to be {\em finitely generated}.
\fresult

We remark that, in the sum above, each $\alpha_j$ decomposes as a finite $k$-linear combination
$\alpha_j=\sum_{S\in\CC(X,E_j)} \lambda_S\, S$ where $\lambda_S\in k$.
Therefore, every $m\in M(X)$ decomposes further as a (finite) $k$-linear combination
$$m=\sumb{j\in J}{\rule{0ex}{1.6ex}S\in\CC(X,E_j)} \lambda_S\, S m_j \mpoint$$

\result{Example} \label{Ex-rep}
{\rm If $E$ is a finite set, the representable functor $k\CC(-,E)$ is finitely generated. It is actually generated by a single element, namely $\Delta_E\in k\CC(E,E)$.}
\fresult

\result{Lemma} \label{finite-dim} Let $M$ be a finitely generated correspondence functor over~$k$.
Then, for every finite set~$X$, the evaluation $M(X)$ is a finitely generated $k$-module.
In particular, if $k$ is a field, then $M(X)$ is finite-dimensional.
\fresult

\pf
Let $\{m_i\mid i=1,\ldots,n\}$ be a finite set of generators of~$M$, with $m_i\in M(E_i)$.
Let $B_X= \{Sm_i \mid S\in\CC(X,E_i),  i=1,\ldots,n\}$.
By definition and by the remark above, every element of~$M(X)$ is a $k$-linear combination of elements of~$B_X$.
But $B_X$ is a finite set, so $M(X)$ is finitely generated. If $k$ is a field, this means that $M(X)$ is finite-dimensional.
\endpf

It follows that, in order to understand finitely generated correspondence functors, we could assume that all their evaluations are finitely generated $k$-modules. But we do not need this for our next characterizations.

\result{Proposition} \label{fg} Let $M$ be a correspondence functor over $k$.
The following conditions are equivalent~:
\begin{enumerate}
\item $M$ is finitely generated.
\item $M$ is isomorphic to a quotient of a functor of the form
$\bigoplusl_{i=1}^n k\CC(-,E_i)$ for some finite sets $E_i$ ($i=1,\ldots,n$).
\item $M$ is isomorphic to a quotient of a functor of the form $\bigoplusl_{i\in I} kC(-,E)$
for some finite set $E$ and some finite index set $I$.
\item There exists a finite set $E$ and a finite subset $B$ of~$M(E)$ such that $M$ is generated by~$B$.
\end{enumerate}
\fresult

\pf
(a) $\Rightarrow$ (b). Suppose that $M$ is generated by the set $\{m_i\mid i=1,\ldots,n\}$, where $m_i\in M(E_i)$.
It follows from Yoneda's lemma that there is a morphism
$$\psi_i:k\CC(-,E_i)\to M$$
mapping $\Delta_{E_i}\in \CC(E_i,E_i)$ to the element $m_i\in M(E_i)$,
hence mapping $\beta\in \CC(X,E_i)$ to $\beta m_i\in M(X)$. Their sum yields a morphism
$$\psi:\bigoplusl_{i=1}^n k\CC(-,E_i)\to M \mpoint$$
For any~$X$ and any $m\in M(X)$,
we have $m=\sum\limits_{i=1}^n \alpha_im_i$ for some $\alpha_i\in k\CC(X,E_i)$,
hence $m=\psi(\alpha_1,\ldots,\alpha_n)$, proving the surjectivity of~$\psi$.\mpn

(b) $\Rightarrow$ (c). Suppose that $M$ is isomorphic to a quotient of a functor of the form
$\bigoplusl_{i=1}^n k\CC(-,E_i)$. Let $F$ be the largest of the sets~$E_i$.
By Corollary~\ref{summand2}, each $k\CC(-,E_i)$ is a direct summand of~$k\CC(-,F)$.
Therefore, $M$ is also isomorphic to a quotient of the functor $\bigoplusl_{i=1}^n k\CC(-,F)$.\mpn

(c) $\Rightarrow$ (d). By Example~\ref{Ex-rep}, $k\CC(-,E)$ is generated by $\Delta_E\in k\CC(E,E)$.
Let $b_i\in \bigoplusl_{i\in I} k\CC(E,E)$ having zero components everywhere, except the $i$-th equal to~$\Delta_E$.
Since $M$ is a quotient of $\bigoplusl_{i\in I} k\CC(-,E)$, it is generated by the images of the elements~$b_i$.
This is a finite set because $I$ is finite by assumption.\mpn

(d) $\Rightarrow$ (a). Since $M$ is generated by~$B$, it is finitely generated.
\endpf

We apply this to the functors $L_{E,V}$ and $S_{E,V}$ defined in Lemma~\ref{adjunction} and Notation~\ref{SXW}.

\result{Corollary} \label{LEV-fg} Let $V$ be a finitely generated $\CR_E$-module,
where $E$ is a finite set and $\CR_E$ is the algebra of relations on~$E$.
Then $L_{E,V}$ and $S_{E,V}$ are finitely generated correspondence functors.
\fresult

\pf
Let $\{v_i\mid i\in I\}$ be a finite set of generators of~$V$ as an $\CR_E$-module.
There is a morphism $\pi_i:k\CC(-,E)\to L_{E,V}$ mapping $\Delta_E$ to $v_i\in L_{E,V}(E)=V$.
Therefore, we obtain a surjective morphism
$$\sum_{i\in I}\pi_i : \bigoplusl_{i\in I} k\CC(-,E)\longrightarrow L_{E,V} \mvirg$$
showing that $L_{E,V}$ is finitely generated.
Now $S_{E,V}$ is a quotient of $L_{E,V}$, so it is also finitely generated.
\endpf

\result{Proposition}\label{KRS} Let $k$ be a noetherian ring.
\begin{enumerate} 
\item For any finitely generated correspondence functor $M$ over $k$, the algebra $\End_{\CF_k}(M)$ is a finitely generated $k$-module.
\item For any two finitely generated correspondence functors $M$ and $N$, the $k$-module $\Hom_{\CF_k}(M,N)$ is finitely generated.
\item If $k$ is a field, the Krull-Remak-Schmidt theorem holds for finitely generated correspondence functors over~$k$. 
\end{enumerate}
\fresult

\pf (a) Since $M$ is finitely generated, there exists a finite set $E$ and a surjective morphism
$\pi:\bigoplusl_{i\in I} k\mathcal{C}(-,E)\to M$ for some finite set~$I$ (Proposition~\ref{fg}). Denote by $\mathcal{A}$ the subalgebra of $\End_{\mathcal{F}_k}\big(\bigoplusl_{i\in I} k\mathcal{C}(-,E)\big)$ consisting of endomorphisms $\varphi$ such that $\varphi(\Ker\pi)\subseteq\Ker\pi$. The algebra $\mathcal{A}$ is isomorphic to a $k$-submodule
of $\End_{\mathcal{F}_k}\big(\bigoplusl_{i\in I} k\mathcal{C}(-,E)\big)$, which is isomorphic to a matrix algebra of size~$|I|$ over the $k$-algebra $k\mathcal{C}(E,E)$ (because $\End_{\CF_k}\big(k\mathcal{C}(-,E)\big)\cong k\mathcal{C}(E,E)$ by Yoneda's lemma). This matrix algebra is free of finite rank as a $k$-module. As $k$ is noetherian, it follows that $\mathcal{A}$ is a finitely generated $k$-module.\par

Now by definition of~$\mathcal{A}$, any $\varphi\in \mathcal{A}$ induces an endomorphism $\sur{\varphi}$ of $M$ such that $\sur{\varphi}\pi=\nolinebreak\pi\varphi$. This yields an algebra homomorphism $\mathcal{A}\to\End_{\mathcal{F}_k}(M)$, which is surjective, since the functor $k\mathcal{C}(-,E)$ is projective. It follows that $\End_{\mathcal{F}_k}(M)$ is also a finitely generated $k$-module.\mpn

(b) The functor $M\oplus N$ is finitely generated, hence $V=\End_{\CF_k}(M\oplus N)$ is a finitely generated $k$-module, by~(a). Since $\Hom_{\CF_k}(M,N)$ embeds in $V$, it is also a finitely generated $k$-module.\mpn

(c) If moreover $k$ is a field, then $\End_{\mathcal{F}_k}(M)$ is a finite dimensional $k$-vector space, by~(a). Finding decompositions of $M$ as a direct sum of subfunctors amounts to splitting the identity of $\End_{\mathcal{F}_k}(M)$ as a sum of orthogonal idempotents. Since $\End_{\mathcal{F}_k}(M)$ is a finite dimensional algebra over the field $k$, the standard theorems on decomposition of the identity as a sum of primitive idempotents apply.
Thus $M$ can be split as a direct sum of indecomposable functors, and such a decomposition is unique up to isomorphism.
\endpf

After the Krull-Remak-Schmidt theorem, we treat the case of projective covers.
Recall (see 2.5.14 in~\cite{AF} for categories of modules) that in an abelian category~$\mathcal{A}$, a subobject $N$ of an object $P$ is called {\em superfluous} if for any subobject $X$ of $P$, the equality $X+N=P$ implies $X=P$. Similarly, an epimorphism $f:P\to M$ in $\mathcal{A}$ is called superfluous if $\Ker f$ is superfluous in $P$, or equivalently, if for any morphism $g:L\to P$ in $\mathcal{A}$, the composition $f\circ g$ is an epimorphism if and only if $g$ is an epimorphism. A {\em projective cover} of an object $M$ of $\mathcal{A}$ is defined as a pair $(P,p)$, where $P$ is projective and $p$ is a superfluous epimorphism from $P$ to~$M$.

\result{Proposition} \label{proj-cover}
Let $M$ be a finitely generated correspondence functor over a commutative ring~$k$.
\begin{enumerate}
\item Suppose that $M$ is generated by $M(E)$ where $E$ is a finite. If $(P,p)$ is a projective cover of $M(E)$ in $\CR_E\text{-\!}\Mod$, then $(L_{E,P},\widetilde{p})$ is a projective cover of $M$ in $\mathcal{F}_k$, where $\widetilde{p}:L_{E,P}\to M$ is obtained from $p:P\to M(E)$ by the adjunction of Lemma~\ref{adjunction}.
\item If $k$ is a field, then $M$ admits a projective cover.
\item In particular, when $k$ is a field, let $E$ be a finite set, let $R$ be an order relation on $E$, and let $V$ be a simple $k\Aut(E,R)$-module. Let moreover $(P,p)$ be a projective cover of $\CP f_R\otimes_{k\Aut(E,R)}V$ in $\CR_E\text{-\!}\Mod$. Then $(L_{E,P},\widetilde{p})$ is a projective cover of the simple correspondence functor $S_{E,R,V}$.
\end{enumerate}
\fresult

\pf (a) (This was already proved in Lemme 2 of~\cite{Bo1}.) By Lemma~\ref{adjunction}, the functor $Q\mapsto L_{E,Q}$ maps projectives to projectives. So the functor $L_{E,P}$ is projective. Since $M$ is generated by $M(E)$, and since the evaluation at $E$ of the morphism $\widetilde{p}:L_{E,P}\to M$ is equal to $p:P\to M(E)$, it follows that $\widetilde{p}$ is surjective. If $N$ is any subfunctor of $L_{E,P}$ such that $\widetilde{p}(N)=M$, then in particular $N(E)\subseteq P$ and $p\big(N(E)\big)=M(E)$. Since $p$ is superfluous, it follows that $N(E)=P$, hence $N=L_{E,P}$ since $L_{E,P}$ is generated by its evaluation $P$ at $E$.\mpn

(b) The algebra $\CR_E$ is a finite dimensional algebra over the field~$k$. Hence any finite dimensional $\CR_E$-module admits a projective cover. Therefore (b) follows from~(a).\mpn

(c) The evaluation of the simple functor $S_{E,R,V}$ at $E$ is the simple $\CR_E$-module $\CP f_R\otimes_{k\Aut(E,R)}V$. Hence (c) follows from (a) and (b).
\endpf


\section{Bounded type} \label{Section-bounded}

\bigskip
\noindent
In this section, we analyze a notion which is more general than finite generation.

\result{Definition} Let $k$ be a commutative ring and let $M$ be a correspondence functor over $k$. 
\begin{enumerate}
\item We say that $M$ has {\em bounded type} if there exists a finite set $E$ such that $M$ is generated by $M(E)$.
\item We say that $M$ has a {\em bounded presentation} if there are projective correspondence functors $P$ and $Q$ of bounded type and an exact sequence of functors
$$Q\to P\to M\to 0 \mpoint$$
Such a sequence is called a {\em bounded presentation} of $M$.
\end{enumerate}
\fresult

Suppose that $M$ has bounded type and let $E$ be a finite set such that $M$ is generated by~$M(E)$.
It is elementary to see that $M$ is finitely generated if and only if $M(E)$ is a finitely generated $\CR_E$-module (using Example~\ref{Ex-rep} and Lemma~\ref{finite-dim}). Thus an infinite direct sum of copies of a simple functor $S_{E,R,V}$ has bounded type (because it generated by its evaluation at~$E$) but is not finitely generated.
Also, a typical example of a correspondence functor which does not have bounded type is a direct sum of simple functors
$\bigoplusl_{n=0}^\infty S_{E_n,R_n,V_n}$, where $|E_n|=n$ for each~$n$. This is because $S_{E_n,R_n,V_n}$ cannot be generated by a set of cardinality~$<n$.

\result{Lemma} \label{F>E} Let $k$ be a commutative ring and let $M$ be a correspondence functor over~$k$.
Suppose that $M$ has bounded type and let $E$ be a finite set such that $M$ is generated by~$M(E)$.
For any finite set $F$ with $|F|\geq |E|$, the functor $M$ is generated by $M(F)$.
\fresult

\pf 
Let $i_*\in\CC(F,E)$ and $i^*\in \CC(E,F)$ be as in Lemma~\ref{summand1},
so that $i^*i_*=\Id_E$.
Saying that $M$ is generated by $M(E)$ amounts to saying that $M(X)$ is equal to $k\CC(X,E)M(E)$, for any finite set $X$. It follows that
\begin{eqnarray*}
M(X)&=&k\mathcal{C}(X,E)M(E)\\
&=&k\mathcal{C}(X,E)i^*i_*M(E)\\
&\subseteq& k\mathcal{C}(X,F)i_*M(E)\\
&\subseteq& k\mathcal{C}(X,F)M(F) \subseteq M(X)\mvirg
\end{eqnarray*}
hence $M(X)=k\mathcal{C}(X,F)M(F)$, i.e. $M$ is generated by $M(F)$.
\endpf

We are going to prove that any correspondence functor having a bounded presentation is isomorphic to some functor $L_{E,V}$.
We first deal with the case of projective functors.

\result{Lemma} \label{counit-for-proj} 
Suppose that a correspondence functor $M$ has bounded type and let $E$ be a finite set such that $M$ is generated by~$M(E)$.
If $M$ is projective, then for any finite set $F$ with $|F|\geq |E|$, the $\CR_F$-module $M(F)$ is projective, and the counit morphism $L_{F,M(F)}\to M$ is an isomorphism.
\fresult

\pf
By Lemma~\ref{F>E}, $M$ is generated by $M(F)$.
Choosing a set $B$ of generators of $M(F)$ as an $\CR_F$-module (e.g.~$B=M(F)$), we see that $M$ is also generated by~$B$.
As in the beginning of the proof of Proposition~\ref{fg}, we can apply Yoneda's lemma and obtain a surjective morphism
$\bigoplusl_{b\in B} k\CC(-,F)\to M$. Since $M$ is projective, this morphism splits, and its evaluation at~$F$ also splits as a map of $\CR_F$-modules. Hence $M(F)$ is isomorphic to a direct summand of a free $\CR_F$-module, that is, a  projective $\CR_F$-module.\par

By adjunction (Lemma~\ref{adjunction}), there is a morphism $\theta:L_{F,M(F)}\to M$ which, evaluated at $F$, gives the identity map of~$M(F)$. As $M$ is generated by~$M(F)$, it follows that $\theta$ is surjective, hence split since $M$ is projective. Let $\eta:M\to L_{F,M(F)}$ be a section of~$\theta$. Since, on evaluation at~$F$, we have $\theta_F=\Id_{M(F)}$, the equation $\theta\eta=\Id_M$ implies that, on evaluation at~$F$, we get $\eta_F=\Id_{M(F)}$.
Therefore $\eta_F\theta_F=\Id_{M(F)}$.
Now $\eta\theta: L_{F,M(F)}\to L_{F,M(F)}$ corresponds by adjunction to
$$\Id_{M(F)} = \eta_F\theta_F: M(F) \longrightarrow L_{F,M(F)}(F) = M(F) \mpoint$$
Therefore $\eta\theta$ must be the identity. It follows that $\eta$ and $\theta$ are mutual inverses.
Thus $M\cong L_{F,M(F)}$.
\endpf

We now prove that any functor with a bounded presentation is an $L_{E,V}$, and conversely.
In the case of a noetherian base ring~$k$, this result will be improved in Section~\ref{Section-noetherian}.

\result{Theorem} \label{counit-iso}
\begin{enumerate}
\item Suppose that a correspondence functor $M$ has a bounded presentation
$$Q\to P\to M\to 0 \mpoint$$
Let $E$ be a finite set such that $P$ is generated by $P(E)$ and $Q$ is generated by $Q(E)$.
Then for any finite set $F$ with $|F|\geq |E|$, the counit morphism
$$\eta_{M,F}:L_{F,M(F)}\to M$$
is an isomorphism.
\item If $E$ is a finite set and $V$ is an $\CR_E$-module, then
the functor $L_{E,V}$ has a bounded presentation. More precisely, if
$$W_1\to W_0\to V\to 0$$
is a projective resolution of $V$ as an $\CR_E$-module, then
$$L_{E,W_1}\to L_{E,W_0}\to L_{E,V}\to 0$$
is a bounded presentation of $L_{E,V}$.
\end{enumerate}
\fresult

\pf 
(a) Consider the commutative diagram
$$\xymatrix{
L_{F,Q(F)}\ar[d]^-{\eta_{Q,F}}\ar[r]&L_{F,P(F)}\ar[d]^-{\eta_{P,F}}\ar[r]&L_{F,M(F)}\ar[d]^-{\eta_{M,F}}\ar[r]&0\\
Q\ar[r]&P\ar[r]&M\ar[r]&0
}$$
where the vertical maps are obtained by the adjunction of Lemma~\ref{adjunction}.
This lemma also asserts that the first row is exact.
By Lemma~\ref{counit-for-proj}, for any finite set $F$ with $|F|\geq |E|$,
the vertical morphisms $\eta_{Q,F}$ and $\eta_{P,F}$ are isomorphisms.
Since the rows of this diagram are exact, it follows that $\eta_{M,F}$ is also an isomorphism.\mpn

(b) We use the adjunction of Lemma~\ref{adjunction}. Applying the right exact functor $U\mapsto L_{E,U}$ to the exact sequence $W_1\to W_0\to V\to 0$ gives the exact sequence
$$L_{E,W_1}\to L_{E,W_0}\to L_{E,V}\to 0\mpoint$$
By Lemma~\ref{adjunction}, $L_{E,W_1}$ and $L_{E,W_0}$ are projective functors, since $W_1$ and $W_0$ are projective $\CR_E$-modules. They all have bounded type since they are generated by their evaluation at~$E$.
\endpf

Given a finite set $E$ and an $\CR_E$-module $V$, we define an induction procedure as follows.
For any finite set~$F$, we define the $\CR_F$-module
$$V\!\uparrow_E^F \;:=k\CC(F,E)\otimes_{\CR_E}V \mpoint$$
Notice that, by the definition of~$L_{E,V}$, we have $L_{E,V}(F)=V\!\uparrow_E^F$.
To end this section, we mention the behavior of the functors $L_{E,V}$ under induction.

\result{Proposition} \label{induction} Let $E$ be a finite set and $V$ be an $\CR_E$-module. 
If $F$ is a finite set with $|F|\geq |E|$, the equality $L_{E,V}(F)= {V\!\uparrow_E^F}$ induces an isomorphism of correspondence functors
$$L_{F,V\!\uparrow_E^F}\cong L_{E,V}\mpoint$$
\fresult

\pf Let $M=L_{E,V}$. Then by Theorem~\ref{counit-iso}, there exists a bounded presentation
$$Q\to P\to M\to 0$$
where $Q=L_{E,W_1}$ is generated by $Q(E)$ and $P=L_{E,W_0}$ is generated by $P(E)$.
Hence by Theorem~\ref{counit-iso}, for any finite set $F$ with $|F|\geq |E|$, the counit morphism 
$$\eta_{M,F}:L_{F,M(F)}\to M$$
is an isomorphism. In other words $\eta_{M,F}:L_{F,V\!\uparrow_E^F}\to L_{E,V}$ is an isomorphism.\endpf


\section{Exponential behavior and finite length} \label{Section-exponential}

\bigskip
\noindent
In this section, we give a lower bound estimate for the dimension of the evaluations of a simple functor $S_{E,R,V}$,
which is proved to behave exponentially.
We also prove that the exponential behavior is equivalent to finite generation. 

We first need a well-known combinatorial lemma.

\result{Lemma} \label{sandwich} Let $E$ be a finite set and let $G$ be a finite set containing~$E$.
\begin{enumerate}
\item For any finite set $X$, the number $s(X,E)$ of surjective maps $\varphi:X\to E$ is equal to
$$s(X,E)=\sum_{i=0}^{|E|}(-1)^{|E|-i}{|E|\choose i} \;i^{|X|} \mvirg$$
or equivalently
$$s(X,E)=\sum_{j=0}^{|E|}(-1)^j{|E|\choose j} \;(|E|-j)^{|X|} \mpoint$$
\item More generally, for any finite set $X$, the number $ss(X,E,G)$ of all maps ${\varphi:X\to G}$ such that $E\subseteq \varphi(X)\subseteq G$ is equal to
$$ss(X,E,G)=\sum_{i=0}^{|E|}(-1)^i{|E|\choose i}(|G|-i)^{|X|}\mpoint$$
\end{enumerate}
\fresult

\pf (a) Up to multiplication by $|E|!$, the number $s(X,E)$ is known as a Stirling number of the second kind.   
Either by Formula~(24a) in Section~1.4 of~\cite{St}, or by a direct application of M\"obius inversion (i.e.~inclusion-exclusion principle in the present case), we have
$$s(X,E)=\sum_{B\subseteq E}(-1)^{|E-B|}|B|^{|X|}\mpoint$$
Setting $|B|=i$, the first formula in (a) follows.\mpn

(b) Applying (a) to each subset $A$ such that $E\subseteq A\subseteq G$, we obtain
\begin{eqnarray*}
ss(X,E,G)&=&\sum_{E\subseteq A\subseteq G}s(X,A)\\
&=&\sum_{E\subseteq A\subseteq G}\sum_{B\subseteq A}(-1)^{|A-B|}|B|^{|X|}\\
&=&\sum_{B\subseteq G}\Big(\sum_{E\cup B\subseteq A\subseteq G}(-1)^{|A-B|}\Big)|B|^{|X|}
\end{eqnarray*}
But the inner sum is zero unless $E\cup B=G$. Therefore
\begin{eqnarray*}
ss(X,E,G)&=&\sumb{B\subseteq G}{\rule{0ex}{1.6ex}E\cup B=G}(-1)^{|G-B|}|B|^{|X|}\\
&=&\sum_{C\subseteq E}(-1)^{|C|}\big(|G-C|\big)^{|X|}\mvirg
\end{eqnarray*}
where the last equality follows by setting $C=G-B$. This proves part~(b).
\endpf

Now we prove our main lower bound estimate for the dimensions of the evaluations of a simple functor.

\result{Theorem} \label{exponential-simple}
Suppose that $k$ is a field and let $S_{E,R,V}$ be a simple correspondence functor, where $E$ is a finite set,
$R$ is an order on~$E$, and $V$ is a simple $k\Aut(E,R)$-module.
There exists a positive integer~$N$ and a positive real number~$c$ such that, for any finite set~$X$ of cardinality at least~$N$, we have
$$c |E|^{|X|} \leq \dim\big(S_{E,R,V}(X)\big) \leq \big(2^{|E|}\big)^{|X|} \mpoint$$
\fresult

\pf
Recall that $S_{E,R,V}=L_{E,T_{R,V}}/J_{E,T_{R,V}}$ where $T_{R,V}$ denotes the $\CR_E$-module
$$T_{R,V}=\CP_Ef_R\otimes_{k\Aut(E,R)}V \mpoint$$
Since $V$ is simple, it is generated by a single element $v\in V-\{0\}$.
Since $f_R\alpha=\alpha f_R$ for any $\alpha\in\Aut(E,R)$, it follows that
the $\CR_E$-module $\CP_Ef_R\otimes_{k\Aut(E,R)}V$ is generated by the single element $f_R\otimes v$.
Therefore, we have a surjective morphism of correspondence functors
$$\begin{array}{rcl}
k\CC(-,E) &\longrightarrow&
k\CC(-,E) \otimes_{\CR_E} \CP_Ef_R\otimes_{k\Aut(E,R)}V
=  L_{E,T_{R,V}} \\
U & \mapsto & U\otimes f_R\otimes v \mpoint
\end{array}$$
Since $S_{E,R,V}$ is a quotient of $L_{E,T_{R,V}}$, we obtain a surjective morphism of correspondence functors
$$k\CC(-,E) \longrightarrow S_{E,R,V} \,, \qquad U \mapsto \overline{U\otimes f_R\otimes v} \mvirg$$
where, for any $U\in\CC(X,E)$, the element $\overline{U\otimes f_R\otimes v}\in S_{E,R,V}(X)$ denotes the class of
$U\otimes f_R\otimes v\in L_{E,T_{R,V}}(X)$.

We first prove the upper bound. This is easy and holds for every finite set~$X$.
Since $S_{E,R,V}$ is isomorphic to a quotient of $k\CC(-,E)$, we have
$$\dim\big(S_{E,R,V}(X)\big) \leq \dim\big(k\CC(X,E)\big) = 2^{|X\times E|}  = \big(2^{|E|}\big)^{|X|} \mpoint$$

In order to find a lower bound, for some finite set~$X$,
we introduce the set $\Phi$ of all surjective maps $\varphi: X\to E$.
The symmetric group $\Sigma_E$ acts (on the left) on~$\Phi$ by composition.
Since $\Phi$ consists of surjections onto~$E$, this action is free, that is, the stabilizer of each $\varphi\in\Phi$ is trivial.
We consider the subgroup $\Aut(E,R)$ of $\Sigma_E$
and we let $A$ be a set of representatives of the set of left orbits $\Aut(E,R)\dom \Phi$.\par

For any $\varphi\in\Phi$, we define
$$\Lambda_\varphi = \{ (x,e)\in X\times E \mid (\varphi(x),e)\in R \} \quad \text{and} \quad
\Gamma_\varphi = \{ (e,x)\in E\times X \mid (e,\varphi(x))\in R \}$$
We claim that $\Lambda_\varphi R=\Lambda_\varphi$ and $R\Gamma_\varphi=\Gamma_\varphi$.
Since $\Delta_E\subseteq R$, we always have
$\Lambda_\varphi=\Lambda_\varphi\Delta_E\subseteq\Lambda_\varphi R$.
Conversely, if $(x,f)\in \Lambda_\varphi R$, then there exists $e\in E$ such that
$(x,e)\in \Lambda_\varphi$ and $(e,f)\in R$, that is, $(\varphi(x),e)\in R$ and $(e,f)\in R$.
It follows that $(\varphi(x),f)\in R$ by transitivity of~$R$, that is, $(x,f)\in \Lambda_\varphi$.
Thus $\Lambda_\varphi R\subseteq \Lambda_\varphi$ and equality follows.
The proof for $\Gamma_\varphi$ is similar.\par

Now we consider the set
$$\{ \Lambda_\varphi \mid \varphi\in A \} \subseteq \CC(X,E) \mpoint$$
We want to prove that the image of this set in $S_{E,R,V}(X)$ is linearly independent,
from which we will deduce that $|A|\leq \dim\big(S_{E,R,V}(X)\big)$.
Suppose that $\sum\limits_{\varphi\in A}\lambda_\varphi \Lambda_\varphi$ is mapped to zero in $S_{E,R,V}(X)$,
where $\lambda_\varphi\in k$ for every $\varphi\in A$. In other words,
$$\sum\limits_{\varphi\in A}\lambda_\varphi \overline{\Lambda_\varphi\otimes f_R\otimes v}=0 \mvirg
\qquad \text{that is,} \qquad
\sum\limits_{\varphi\in A}\lambda_\varphi \Lambda_\varphi\otimes f_R\otimes v \in J_{E,T_{R,V}} \mpoint$$
The definition of~$J_{E,T_{R,V}}$ implies that, for every $U\in\CC(E,X)$,
$$\sum\limits_{\varphi\in A}\lambda_\varphi U\Lambda_\varphi\cdot( f_R\otimes v) =0 \mpoint$$
Choosing in particular $U=\Gamma_\psi$ and $\psi\in A$, we obtain~:
\begin{equation} \label{action-Gamma-Lambda}
\text{for every } \; \psi\in A \,, \qquad
\sum\limits_{\varphi\in A}\lambda_\varphi \Gamma_\psi\Lambda_\varphi\cdot( f_R\otimes v) =0 \mpoint
\end{equation}
By Proposition~\ref{fundamental-module}, the action of the relation $\Gamma_\psi\Lambda_\varphi$ on~$f_R\otimes v$ is given by
$$(\Gamma_\psi\Lambda_\varphi)\cdot f_R\otimes v
=\left\{\begin{array}{ll}
\Delta_{\tau}f_R\otimes v&\hbox{if}\;\;\exists\tau\in\Sigma_E\;\hbox{such that}\;
\Delta_E\subseteq \Delta_{\tau^{-1}}\Gamma_\psi\Lambda_\varphi\subseteq R \,,\\
0&\hbox{otherwise}\,.\end{array}\right.$$
In the first case, $\tau$ is unique.\par

We claim that
$$\Delta_E\subseteq \Delta_{\tau^{-1}}\Gamma_\psi\Lambda_\varphi\subseteq R
\; \iff \; \Gamma_\psi\Lambda_\varphi= \Delta_\tau R \mpoint$$
If the left hand side holds, then multiply on the right by~$R$ and use the fact that $\Lambda_\varphi R=\Lambda_\varphi$
to obtain $\Delta_{\tau^{-1}}\Gamma_\psi\Lambda_\varphi= R$, hence $\Gamma_\psi\Lambda_\varphi= \Delta_\tau R$.
Conversely, if the right hand side holds, then $\Delta_{\tau^{-1}}\Gamma_\psi\Lambda_\varphi = R$, hence
$$R\Delta_{\tau^{-1}}\Gamma_\psi\Lambda_\varphi= R^2=R$$
by transitivity and reflexivity of~$R$. In particular, by reflexivity again,
$$\Delta_E\subseteq R\Delta_{\tau^{-1}}\Gamma_\psi\Lambda_\varphi$$
so that, for any $(a,a)\in\Delta_E$, there exists $b\in E$ with $(a,b)\in R$
and $(b,a)\in\Delta_{\tau^{-1}}\Gamma_\psi\Lambda_\varphi=R$.
By antisymmetry of~$R$, it follows that $b=a$ and therefore
$(a,a)\in\Delta_{\tau^{-1}}\Gamma_\psi\Lambda_\varphi$,
so that $\Delta_E\subseteq \Delta_{\tau^{-1}}\Gamma_\psi\Lambda_\varphi$.
This shows that the left hand side holds, proving the claim.\par

We can now rewrite (\ref{action-Gamma-Lambda}) as follows~:
$$\text{for every } \; \psi\in A \,, \qquad
\sum_{\substack{\varphi\in A\,,\,\tau\in\Sigma_E \\
\Gamma_\psi\Lambda_\varphi= \Delta_\tau R}}
\lambda_\varphi \Delta_\tau f_R\otimes v =0 \mpoint$$
Write $\tau=\sigma\alpha$ with $\sigma\in S$ and $\alpha\in\Aut(E,R)$, where $S$ denotes some set of representatives of cosets $\Sigma_E/\Aut(E,R)$. Since $\Delta_\alpha f_R=f_R\Delta_\alpha$ for every $\alpha\in\Aut(E,R)$ and since  the tensor product is over~$k\Aut(E,R)$, we obtain
\begin{equation} \label{equation-for-all-psi}
\text{for every } \; \psi\in A \,, \qquad
\sum_{\sigma\in S}  \Delta_\sigma f_R\otimes \Big(
\sum_{\substack{\varphi\in A\,,\,\alpha\in\Aut(E,R) \\
\Gamma_\psi\Lambda_\varphi= \Delta_{\sigma\alpha} R}}
\lambda_\varphi \Delta_\alpha v\Big) =0 \mpoint
\end{equation}
By Proposition~\ref{fundamental-module}, $\CP_E f_R$ is a $(\CP_E,k\Aut(E,R))$-bimodule with a $k$-basis consisting of the elements $\Delta_\tau f_R$, for $\tau\in\Sigma_E$. Therefore,
$$\CP_E f_R =\bigoplus_{\sigma\in S} \Delta_\sigma f_R \cdot(k\Aut(E,R))$$
and it follows that
$$\CP_E f_R \otimes_{k\Aut(E,R)} V =\bigoplus_{\sigma\in S} \Delta_\sigma f_R \otimes_{k\Aut(E,R)} V \mpoint$$
Since the sum is direct and since $f_R \otimes_{k\Aut(E,R)} V\cong V$ because the right action of~$\Aut(E,R)$ on~$f_R$ is free, each inner sum in (\ref{equation-for-all-psi}) is zero. In particular, taking $\sigma=\Id$ (which we may choose in~$S$), we get~:
\begin{equation} \label{equation-for-all-psi-2}
\text{for every } \; \psi\in A \,, \qquad
\sum_{\substack{\varphi\in A\,,\,\alpha\in\Aut(E,R) \\
\Gamma_\psi\Lambda_\varphi= \Delta_\alpha R}}
\lambda_\varphi \Delta_\alpha v =0 \mpoint
\end{equation}

Let us analyze the condition $\Gamma_\psi\Lambda_\varphi= \Delta_\alpha R$.
For any given $x\in X$, we can choose $e=\psi(x)$ and $f=\varphi(x)$
and we obtain $(e,x)\in \Gamma_\psi$ and $(x,f)\in \Lambda_\varphi$, hence $(e,f)\in \Gamma_\psi\Lambda_\varphi=\Delta_\alpha R$, that is, $(\alpha^{-1}(e),f)\in R$.
In other words, if we write simply $\leq_R$ for the partial order~$R$, we obtain $\alpha^{-1}(\psi(x))\leq_R \varphi(x)$.
This holds for every $x\in X$ and we define
$$\psi\preceq\varphi \;\iff\; \exists\, \alpha\in\Aut(E,R) \; \text{ such that }\;
\alpha^{-1}(\psi(x))\leq_R \varphi(x)\quad \text{for all } \,x\in X \mpoint$$
Therefore the condition $\Gamma_\psi\Lambda_\varphi= \Delta_\alpha R$ implies that $\psi\preceq\varphi$.\par

Now we prove that the relation $\preceq$ is a partial order on the set~$A$. It is reflexive, by taking simply $\alpha=\Id$.
It is transitive because if $\alpha^{-1}(\psi(x))\leq_R \varphi(x)$ and $\beta^{-1}(\varphi(x))\leq_R \chi(x)$ for all~$x$, then
$$(\alpha\beta)^{-1}(\psi(x))=\beta^{-1}\alpha^{-1}(\psi(x))\leq_R \beta^{-1}(\varphi(x))\leq_R\chi(x)$$
using the fact that $\beta^{-1}\in\Aut(E,R)$. Finally, the relation $\preceq$ is antisymmetric because if
$\alpha^{-1}(\psi(x))\leq_R \varphi(x)$ and $\beta^{-1}(\varphi(x))\leq_R \psi(x)$, then
$$(\alpha\beta)^{-1}(\psi(x)) \leq_R \psi(x)$$
from which it follows that
$$\psi(x)=(\alpha\beta)^{-n}(\psi(x)) \leq_R (\alpha\beta)^{-1}(\psi(x)) \leq_R \psi(x) \mvirg$$
where $n$ is the order of~$\alpha\beta$ in the group~$\Aut(E,R)$.
But this implies that $\psi(x)=(\alpha\beta)^{-n}(\psi(x))=(\alpha\beta)^{-1}(\psi(x))$, hence
$$\psi(x)=\beta^{-1}\alpha^{-1}(\psi(x))\leq_R \beta^{-1}(\varphi(x)) \leq_R \psi(x) \mvirg$$
and therefore $\beta^{-1}(\varphi(x)) = \psi(x)$. Thus $\varphi$ and $\psi$ belong to the same orbit under the action of~$\Aut(E,R)$. This forces $\varphi=\psi$ because $\varphi$ and $\psi$ belong to our chosen set~$A$ of representatives of the set of left orbits $\Aut(E,R)\dom \Phi$.\par

In view of proving the linear independence we are looking for, suppose that the coefficients $\lambda_\varphi$ are not all zero. Choose $\psi\in A$ maximal (with respect to~$\preceq$) such that $\lambda_\psi\neq0$.
In the sum (\ref{equation-for-all-psi-2}), the condition $\Gamma_\psi\Lambda_\varphi= \Delta_\alpha R$
implies, as we have seen above, that $\psi\preceq\varphi$. Since $\psi$ is maximal,
the sum over $\varphi\in A$ actually runs over the single element~$\psi$ and reduces to
$$\sum_{\substack{\alpha\in\Aut(E,R) \\
\Gamma_\psi\Lambda_\psi= \Delta_\alpha R}}
\lambda_\psi \Delta_\alpha v =0 \mpoint$$
But the condition $\Gamma_\psi\Lambda_\psi= \Delta_\alpha R$ implies that $\alpha=\Id$ because
$\Delta_E\subseteq \Gamma_\psi\Lambda_\psi$, hence
$\Delta_{\alpha^{-1}} \subseteq \Delta_{\alpha^{-1}}\Gamma_\psi\Lambda_\psi$,
that is, $\Delta_{\alpha^{-1}} \subseteq R$,
which can only occur if $\alpha=\Id$.
It follows that there is a single term in the whole sum (\ref{equation-for-all-psi-2}),
namely $\lambda_\psi v =0$. Since $\lambda_\psi\neq0$, we obtain $v=0$, which is impossible since $v$ was chosen nonzero in~$V$. This contradiction shows that all coefficients $\lambda_\varphi$ are zero,
proving the linear independence of the image of~$A$ in $S_{E,R,V}(X)$.
Therefore $|A|\leq \dim\big(S_{E,R,V}(X)\big)$.\par

Now we need to estimate $|A|$ and, for simplicity, we write $e=|E|$ and $x=|X|$.
Since $A$ is a set of representatives of orbits in~$\Phi$ under the free action of~$\Aut(E,R)$,
we have $|\Phi|=|\Aut(E,R)|\cdot |A|$, so we need to estimate $|\Phi|$.
By Lemma~\ref{sandwich}, we have
$$|\Phi| =\sum_{i=0}^e(-1)^{e-i}{e\choose i} \;i^x
=e^x + \sum_{i=0}^{e-1}(-1)^{e-i}{e\choose i} \;i^x \mpoint$$
Note that the second sum is negative because the number $|\Phi|$ of surjective maps $X\to E$ is smaller than the number $e^x$ of all maps $X\to E$. We can rewrite
$$|\Phi| =e^x \Big(1+\sum_{i=0}^{e-1}(-1)^{e-i}{e\choose i} \;\big(\frac{i}{e}\big)^x\Big) \mpoint$$
Since $\frac{i}{e}\leq \frac{e-1}{e}<1$, the sum can be made as small as we want, provided $x$ is large enough.
Therefore there exists a positive integer~$N$ and a positive real number~$a$ such that $a\,e^x\leq |\Phi|$ whenever $x\geq N$. In other words, for any finite set~$X$ of cardinality at least~$N$, we have
\begin{equation}\label{bound}
\frac a{|\Aut(E,R)|} \, |E|^{|X|} \leq \frac{|\Phi|}{|\Aut(E,R)|} =|A| \leq \dim\big(S_{E,R,V}(X)\big)\mvirg
\end{equation}
giving the required lower bound for $\dim\big(S_{E,R,V}(X)\big)$.
\endpf

We can now characterize finite generation in terms of exponential behavior.

\result{Theorem} \label{exponential}
Let $M$ be a correspondence functor over a field~$k$. The following are equivalent~:
\begin{enumerate}
\item $M$ is finitely generated.
\item There exists positive real numbers $a,b,r$ such that $\dim(M(X))\leq a\, b^{|X|}$ for every finite set~$X$ with $|X|\geq r$.
\end{enumerate}
\fresult

\pf
(a) $\Rightarrow$ (b). Let $M$ be a quotient of $\bigoplusl_{i\in I} k\CC(-,E)$ for some finite set~$E$
and some finite index set~$I$. For every finite set~$X$, we have
$$\dim(M(X))\leq |I| \dim(k\CC(X,E)) = |I| \, 2^{|X\times E|}= |I| \big(2^{|E|}\big)^{|X|} \mpoint$$
\mpn

(b) $\Rightarrow$ (a). Let $P$ and $Q$ be subfunctors of~$M$ such that $Q\subseteq P\subseteq M$ and $P/Q$ simple, hence $P/Q\cong S_{E,R,V}$ for some triple $(E,R,V)$. We claim that $|E|$ is bounded above.
Indeed, for $|X|$ large enough, we have
$$c|E|^{|X|}\leq\dim(S_{E,R,V}(X))$$
for some $c>0$, by Theorem~\ref{exponential-simple}, and
$$\dim(S_{E,R,V}(X)) \leq \dim(M(X))\leq a\, b^{|X|}$$
by assumption. Therefore, whenever $|X|\geq N$ for some~$N$, we have
$$c|E|^{|X|}\leq a\, b^{|X|} \qquad\text{that is, } \quad c\leq a\big(\frac{b}{|E|}\big)^{|X|} \mpoint$$
Since $c>0$, this forces $\displaystyle\frac{b}{|E|}\geq 1$
otherwise $\displaystyle a\big(\frac{b}{|E|}\big)^{|X|}$ is as small as we want.
This shows the bound $|E|\leq b$, proving the claim.\par

For each set $E$ with $|E|\leq b$, we choose a basis $\{m_i\mid 1\leq i\leq n_E\}$ of~$M(E)$
and we use Yoneda's lemma to construct a morphism $\psi_i^E:k\CC(-,E)\to M$ such that, on evaluation at~$E$,
we have $\psi_{i,E}^E(\Delta_E)=m_i$.
Starting from the direct sum of $n_E$ copies of $k\CC(-,E)$, we obtain a morphism
$$\psi^E: k\CC(-,E)^{n_E}\to M \mvirg$$
such that, on evaluation at~$E$, the morphism $\psi_E^E:k\CC(E,E)^{n_E}\to M(E)$
is surjective, because the basis of~$M(E)$ is in the image.
Now the sum of all such morphisms $\psi^E$ yields a morphism
$$\psi: \bigoplusl_{|E|\leq b} k\CC(-,E)^{n_E}\longrightarrow M$$
which is surjective on evaluation at every set~$E$ with $|E|\leq b$.\par

Let $N=\Im(\psi)$ and suppose ab absurdo that $N\neq M$. Let $F$ be a minimal set such that $M(F)/N(F)\neq\{0\}$.
Since $\psi$ is surjective on evaluation at every set~$E$ with $|E|\leq b$, we must have $|F|>b$.
Now $M(F)/N(F)$ is a module for the finite-dimensional algebra $\CR_F=k\CC(F,F)$ and,
by minimality of~$F$, inessential relations act by zero on~$M(F)/N(F)$.
Let $W$ be a simple submodule of $M(F)/N(F)$ as a module for the essential algebra~$\CE_F$.
Associated with~$W$, consider the simple functor $S_{F,W}$.
(Actually, $W$ is parametrized by a pair $(R,V)$ and $S_{F,W}=S_{F,R,V}$ (see Theorem~\ref{simple-modules}), but we do not need this.)
Now the module $W=S_{F,W}(F)$ is isomorphic to a subquotient of $M(F)/N(F)$.
By Proposition~\ref{subquotients}, $S_{F,W}$ is isomorphic to a subquotient of $M/N$.
By the claim proved above, we obtain $|F|\leq b$. This contradiction shows that $N=M$, that is, $\psi$ is surjective.
Therefore $M$ is isomorphic, via~$\psi$, to a quotient of $\bigoplusl\limits_{|E|\leq b} k\CC(-,E)^{n_E}$.
By Proposition~\ref{fg}, $M$ is finitely generated.
\endpf


\section{Finite length} \label{Section-finite-length}

\bigskip
\noindent
Using the exponential behaviour proved in the previous section, we now show that, if our base ring~$k$ is a field, then
every finitely generated correspondence functor has finite length.
We first need a lemma.

\result{Lemma} \label{maximal}
Let $k$ be a field and let $M$ be a finitely generated correspondence functor over~$k$.
\begin{enumerate}
\item $M$ has a maximal subfunctor.
\item Any subfunctor of $M$ is finitely generated.
\end{enumerate}
\fresult

\pf
(a) Since $M$ is finitely generated, $M$ is generated by~$M(E)$ for some finite set~$E$ (Proposition~\ref{fg}).
Let $N$ be a maximal submodule of~$M(E)$ as a $\CR_E$-module. Note that $N$ exists because $M(E)$ is finite-dimensional by Lemma~\ref{finite-dim}. Then $M(E)/N$ is a simple $\CR_E$-module.
By Proposition~\ref{subquotients}, there exist two subfunctors $F\subseteq G\subseteq M$ such that
$G/F$ is simple, $G(E)=M(E)$, and $F(E)=N$.
Since $M$ is generated by~$M(E)$ and $G(E)=M(E)$, we have $G=M$.
Therefore, $F$ is a maximal subfunctor of~$M$.\mpn

(b) Let $N$ be a subfunctor of~$M$. Since $M$ is finitely generated, there exist positive numbers $a,b$ such that,
for every large enough finite set~$X$, we have
$$\dim(N(X))\leq \dim(M(X))\leq a\,b^{|X|} \mvirg$$
by Theorem~\ref{exponential}. The same theorem then implies that $N$ is finitely generated.
\endpf

Lemma~\ref{maximal} fails for other categories of functors.
For instance, in the category of biset functors, the Burnside functor is finitely generated and has a maximal subfunctor which is not finitely generated
(see~\cite{Bo1} or \cite{Bo2}).\par

We now come to one of the most important properties of the category of correspondence functors, namely an artinian property.
As for the previous lemma, the theorem is a specific property of the category of correspondence functors.

\result{Theorem} \label{finite-length}
Let $k$ be a field and let $M$ be a finitely generated correspondence functor over~$k$.
Then $M$ has finite length (that is, $M$ has a finite composition series).
\fresult

\pf
By Lemma~\ref{maximal}, $M$ has a maximal subfunctor~$F_1$ and $F_1$ is again finitely generated.
Then $F_1$ has a maximal subfunctor~$F_2$ and $F_2$ is again finitely generated.
We construct in this way a sequence of subfunctors
\begin{equation}\label{series}
M=F_0\supset F_1\supset F_2 \supset \ldots
\end{equation}
such that $F_i/F_{i+1}$ is simple whenever $F_i\neq 0$. We claim that the sequence is finite,
that is, $F_m=0$ for some~$m$.\par

Let $F_i/F_{i+1}$ be one simple subquotient, hence $F_i/F_{i+1}\cong S_{E,R,V}$ for some triple $(E,R,V)$.
By Theorem~\ref{exponential}, since $M$ is finitely generated, there exist positive numbers $a,b$ such that,
for every large enough finite set~$X$, we have
$$ \dim(M(X))\leq a\,b^{|X|} \mpoint$$
Therefore $\dim(S_{E,R,V}(X))\leq a\,b^{|X|}$.
By Theorem~\ref{exponential-simple}, there exists some constant $c>0$ such that
$c\,|E|^{|X|} \leq \dim(S_{E,R,V}(X))$ for $|X|$ large enough.
So we obtain $c\,|E|^{|X|} \leq a\,b^{|X|}$ for $|X|$ large enough, hence $|E|\leq b$.
This implies that the simple functor $F_i/F_{i+1}\cong S_{E,R,V}$ belongs to a finite set of isomorphism classes of simple functors,
because there are finitely many sets~$|E|$ with $|E|\leq b$ and, for any of them, finitely many order relations~$R$ on~$E$,
and then in turn finitely many $k\Aut(E,R)$-simple modules~$V$ (up to isomorphism).\par

Therefore, if the series~(\ref{series}) of subfunctors~$F_i$ was infinite, then some simple functor $S_{E,R,V}$ would occur infinitely many times (up to isomorphism).
But then, on evaluation at~$E$, the simple $\CR_E$-module $S_{E,R,V}(E)$ would occur
infinitely many times in~$M(E)$. This is impossible because $M(E)$ is finite-dimensional by Lemma~\ref{finite-dim}.
\endpf

Theorem~\ref{finite-length} was obtained independently by Gitlin~\cite{Gi} (for a field $k$ of characteristic zero, or algebraically closed),
using a criterion for finite length proved recently by Wiltshire-Gordon~\cite{WG}.


\section{Projective functors and duality} \label{Section-projective}

\bigskip
\noindent
This section is devoted to projective correspondence functors, mainly in the case where $k$ is a field.
An important ingredient is the use of duality.\par

Recall that, by Lemma~\ref{counit-for-proj}, if a projective correspondence functor $M$ is generated by~$M(E)$, then $M(X)$ is a projective $\CR_X$-module, for every set $X$ with $|X|\geq|E|$.
Recall also that, by Lemma~\ref{counit-for-proj} again, a projective correspondence functor $M$ is isomorphic to $L_{E,M(E)}$ whenever $M$ is generated by~$M(E)$. Thus if we work with functors having bounded type, we can assume that projective functors have the form $L_{E,V}$ for some $\CR_E$-module~$V$. In such a case, we can also enlarge~$E$ because $L_{E,V}\cong L_{F,V\!\uparrow_E^F}$ whenever $|F|\geq|E|$ (see Proposition~\ref{induction}).

\result{Lemma} \label{projective} Let $k$ be a commutative ring and consider the correspondence functor $L_{E,V}$ for some finite set~$E$ and some $\CR_E$-module~$V$.
\begin{enumerate}
\item $L_{E,V}$ is projective if and only if $V$ is a projective $\CR_E$-module.
\item $L_{E,V}$ is finitely generated projective if and only if $V$ is a finitely generated projective $\CR_E$-module.
\item $L_{E,V}$ is indecomposable projective if and only if $V$ is an indecomposable projective $\CR_E$-module.
\end{enumerate}
\fresult

\pf (a) If $L_{E,V}$ is projective, then $V$ is projective by Lemma~\ref{counit-for-proj}. Conversely, if $V$ is projective, then $L_{E,V}$ is projective by Lemma~\ref{adjunction}. \mpn

(b) If $V$ is a finitely generated $\CR_E$-module,
then $L_{E,V}$ is finitely generated by Corollary~\ref{LEV-fg}. If $L_{E,V}$ is finitely generated, then its evaluation $L_{E,V}(E)=V$ is finitely generated by Lemma~\ref{finite-dim}. \mpn

(c) By the adjunction property of Lemma~\ref{adjunction}, $\End_{\mathcal{F}_k}(L_{E,V})\cong\End_{\CR_E}(V)$, so $L_{E,V}$ is indecomposable if and only if $V$ is indecomposable.
\endpf

Our main duality result has two aspects, which we both include in the following theorem.
The notion of symmetric algebra is standard over a field and can be defined over any commutative ring as in~\cite{Br}.

\result{Theorem} \label{self-dual} Let $E$ be a finite set.
\begin{enumerate}
\item The representable functor $k\CC(-,E)$ is isomorphic to its dual.
\item Let $\CR_E=k\CC(E,E)$ be the $k$-algebra of relations on~$E$.
Then $\CR_E$ is a symmetric algebra. More precisely, let $t:\CR_E\to k$ be the $k$-linear form defined, for all basis elements $S \in\CC(E,E)$, by the formula
$$t(S):=\left\{\begin{array}{cl}1&\hbox{if}\;S\cap\Delta_E=\emptyset\mvirg\\
0&\hbox{otherwise}\mpoint\end{array}\right.$$
Then $t$ is a symmetrizing form on~$\CR_E$, in the sense that the associated bilinear form $(a,b)\mapsto t(ab)$ is symmetric and induces an isomorphism of $(\CR_E,\CR_E)$-bimodules between $\CR_E$ and its dual $\Hom_k(\CR_E,k)$.
\end{enumerate}
\fresult

\pf
(a) For every finite set $X$, consider the symmetric bilinear form
$$\big\langle-,-\big\rangle_X : k\CC(X,E) \times k\CC(X,E) \longrightarrow k$$
defined, for all basis elements $R,S \in\CC(X,E)$, by the formula
$$\big\langle R,S \big\rangle_X:=
\left\{ \begin{array}{rl} 1 & \text{if } \, R\cap S=\emptyset \mvirg \\
0 & \text{otherwise.}
\end{array} \right.$$
Then, whenever $U\in \CC(Y,X)$, $R \in\CC(Y,E)$, and $S \in\CC(X,E)$, we have
$$\begin{array}{rl}
R\cap US=\emptyset &\iff \Big( (y,x)\in U, (x,e)\in S \Rightarrow (y,e)\notin R \Big) \\
& \iff \Big( (x,y)\in U\op, (y,e)\in R \Rightarrow (x,e)\notin S \Big) \\
& \iff U\op R \cap S=\emptyset \mpoint
\end{array}$$
It follows that $\big\langle U\op R,S \big\rangle_X=\big\langle R,US \big\rangle_Y$.
In view of the definition of dual functors (Definition~\ref{def-dual}), this implies that the associated family of linear maps
$$\alpha_X:k\CC(X,E) \longrightarrow k\CC(X,E)\dual \,,
\qquad S\mapsto \big(R\mapsto \big\langle R,S \big\rangle_X \big)$$
defines a morphism of correspondence functors $\alpha:k\CC(-,E) \longrightarrow k\CC(-,E)\dual$.\par

To prove that $\alpha$ is an isomorphism, we fix $X$ and we use the complement $\,^c\! R=(X\times E)-R$, for any $R\in \CC(X,E)$.
Notice that the matrix of $\alpha_X$ relative to the canonical basis $\CC(X,E)$ and its dual is the product of two matrices $C$ and $A$, where $C_{R,S}=1$ if $S=\,^c\! R$ and 0 otherwise, while $A$ is the adjacency matrix of the order relation~$\subseteq$. This is because $R\cap S=\emptyset$ if and only if $R\subseteq \,^c\! S$.
Clearly $C$ is invertible (it has order~2) and $A$ is unitriangular, hence invertible. Therefore $\alpha_X$ is an isomorphism.\mpn

(b) Let $R,S\in\CC(E,E)$. Then $t(RS)$ is equal to 1 if $RS\cap\Delta_E=\emptyset$, and $t(RS)=0$ otherwise. Now 
\begin{eqnarray*}
RS\cap\Delta_E=\emptyset&\iff& \Big( (e,f)\in R\ \Rightarrow (f,e)\notin S \Big) \\
&\iff& R\cap S\op=\emptyset \mpoint
\end{eqnarray*}
Therefore $t(RS)=\big\langle R,S\op\big\rangle$,
where $\big\langle-,-\big\rangle_E$ is the bilinear form on $k\CC(E,E)$ defined in~(a).
Since this bilinear form induces an isomorphism with the dual and since the map $S\mapsto S\op$ is an isomorphism (it has order~2), the bilinear form associated with~$t$ induces also an isomorphism with the dual.\par

Since $(R\cap S\op)\op=S\cap R\op$ and $\emptyset\op=\emptyset$, we have $t(RS)=t(SR)$ for any relations $R$ and $S$ on $E$, hence the bilinear form $(a,b)\mapsto t(ab)$ is symmetric.
It is clear that the associated $k$-linear map $\CR_E \to \Hom_k(\CR_E,k)$
is a morphism of $(\CR_E,\CR_E)$-bimodules.
\endpf

\result{Corollary} \label{proj-inj}
If $k$ is a field, then the correspondence functor $k\CC(-,E)$ is both projective and injective.
\fresult

\pf
Since passing to the dual reverses arrows and since $k\CC(-,E)$ is projective, its dual is injective.
But $k\CC(-,E)$ is isomorphic to its dual, so it is both projective and injective.
\endpf

\begin{rem}{Remark} Corollary~\ref{proj-inj} holds more generally when $k$ is a self-injective ring.
\end{rem}

\begin{rem}{Remark} \label{Rdual}
If $R$ is an order relation on~$E$, then there is a direct sum decomposition
$$k\CC(-,E) = k\CC(-,E)R \oplus k\CC(-,E)(1-R) \,.$$
With respect to the bilinear forms defined in the proof of Theorem~\ref{self-dual}, we have
$$\big(k\CC(-,E)R\big)^\perp=k\CC(-,E)(1-R\op)$$
because
$$\begin{array}{rl}
U\in \big(k\CC(X,E)R\big)^\perp &\iff \big\langle U,VR\big\rangle_X=0\;\;\forall V\in k\CC(X,E) \\
&\iff \big\langle UR\op,V\big\rangle_X=0\;\;\forall V\in k\CC(X,E) \\
&\iff UR\op=0 \\
&\iff U(1-R\op)=U \\
&\iff U \in k\CC(X,E)(1-R\op) \mpoint
\end{array}
$$
It follows that
$$k\CC(-,E)/\big(k\CC(-,E)R\big)^\perp =k\CC(-,E)/k\CC(-,E)(1-R\op) \cong k\CC(-,E)R\op$$
and therefore the bilinear forms $\big\langle -,- \big\rangle_X$ induce perfect pairings
$$k\CC(-,E)R \times k\CC(-,E)R\op \longrightarrow k \mpoint$$
Thus $\big(k\CC(-,E)R\big)\dual \cong k\CC(-,E)R\op$.
\end{rem}

\bigskip
By the Krull-Remak-Schmidt theorem (which holds when $k$ is a field by Proposition~\ref{KRS}), it is no harm to assume that our functors are indecomposable.

\result{Theorem} \label{PIF} Let $k$ be a field and $M$ be a finitely generated correspondence functor over $k$.
The following conditions are equivalent:
\begin{enumerate}
\item The functor $M$ is projective and indecomposable.
\item The functor $M$ is projective and admits a unique maximal (proper) subfunctor.
\item The functor $M$ is projective and admits a unique minimal (nonzero) subfunctor.
\item The functor $M$ is injective and indecomposable.
\item The functor $M$ is injective and admits a unique maximal (proper) subfunctor.
\item The functor $M$ is injective and admits a unique minimal (nonzero) subfunctor.
\end{enumerate}
\fresult

\pf
(a) $\Rightarrow$ (b). Suppose first that $M$ is projective and indecomposable. Then $M\cong L_{E,V}$ for some finite set $E$ and some indecomposable projective $\CR_E$-module~$V$ (Lemma~\ref{counit-for-proj}). Since $\CR_E$ is a finite dimensional algebra over~$k$, the module $V$ has a unique maximal submodule~$W$. If $N$ is a subfunctor of $M$, then $N(E)$ is a submodule of $V$, so there are two cases: either $N(E)=V$, and then $N=M$, because $M$ is generated by $M(E)=V$, or $N(E)\subseteq W$, and then $N(X)\subseteq J_W(X)$ for any finite set $X$, where
$$J_W(X)=\Big\{\sum_i\varphi_i\otimes v_i\in k\mathcal{C}(X,E)\otimes_{\CR_E}V\mid\forall \psi\in \CC(E,X),\;\sum_i(\psi\varphi_i)\cdot v_i\in W\Big\}\mpoint$$
One checks easily that the assignment $X\mapsto J_W(X)$ is a subfunctor of $L_{E,V}$, such that $J_W(E)=W$ after the identification $L_{E,V}(E)\cong V$.
(This subfunctor is similar to the one introduced in Lemma~\ref{JXW}.)
In particular $J_W$ is a proper subfunctor of~$L_{E,V}$. It follows that $J_W$ is the unique maximal proper subfunctor of $L_{E,V}$, as it contains any proper subfunctor $N$ of $L_{E,V}$.\mpn

(b) $\Rightarrow$ (a). Suppose that $M$ admits a unique maximal subfunctor $N$. If $M$ splits as a direct sum $M_1\oplus M_2$ of two nonzero subfunctors $M_1$ and $M_2$, then $M_1$ and $M_2$ are finitely generated. Let $N_1$ be a maximal subfunctor of $M_1$, and $N_2$ be a maximal subfunctor of $M_2$. Such subfunctors exist by Lemma~\ref{maximal}. Then $M_1\oplus N_2$ and $N_1\oplus M_2$ are distinct maximal subfunctors of $M$. This contradiction proves that $M$ is indecomposable.\mpn

(a) $\Rightarrow$ (d). If $M$ is a finitely generated projective functor, then there exists a finite set $E$ such that $M$ is isomorphic to a quotient, hence a direct summand, of $\bigoplusl_{i\in I}k\CC(-,E)$ for some finite set~$I$ (Proposition~\ref{fg}). Since $k$ is a field, $k\CC(-,E)$ is an injective functor (Corollary~\ref{proj-inj}), hence so is the direct sum and its direct summand~$M$.\mpn

(d) $\Rightarrow$ (a). If $M$ is a finitely generated injective functor, then its dual $M\dual$ is projective, hence injective, and therefore $M\cong (M\dual)\dual$ is projective.\mpn

(a) $\Rightarrow$ (c). For a finitely generated functor $M$, the duality between $M$ and $M\dual$ induces an order reversing bijection between the subfunctors of $M$ and the subfunctors of~$M\dual$. If $M$ is projective and indecomposable, then so is $M\dual$, that is, (a) holds for~$M\dual$. Thus (b) holds for~$M\dual$ and the functor $M\dual$ has a unique maximal subfunctor. Hence $M$ has a unique minimal subfunctor.\mpn

(c) $\Rightarrow$ (a). If $M$ is projective and admits a unique minimal subfunctor, then $M$ is also injective, and its dual $M\dual$ is projective and admits a unique maximal subfunctor. Hence $M\dual$ is indecomposable, so $M$ is indecomposable.\mpn

It is now clear that (e) and (f) are both equivalent to (a), (b), (c) and~(d).
\endpf

Finally, we prove that the well-known property of indecomposable projective modules over a symmetric algebra also holds for correspondence functors.
Recall that $M/\Rad(M)$ is the largest semi-simple quotient of~$M$ and that $\Soc(M)$ is the largest semi-simple subfunctor of~$M$.

\result{Theorem}\label{dim-Hom-symmetric}
Let $k$ be a field.
\begin{enumerate}
\item Let $M$ be a finitely generated projective correspondence functor over $k$. Then $M/\Rad(M)\cong\Soc(M)$.
\item Let $M$ and $N$ be finitely generated correspondence functors over $k$. If $M$ is projective, then $\dim_k\Hom_{\CF_k}(M,N)=\dim_k\Hom_{\CF_k}(N,M)<+\infty$.
\end{enumerate}
\fresult

\pf (a) By Proposition~\ref{KRS}, we can assume that $M$ is indecomposable. In this case, by Theorem~\ref{PIF}, both $M/\Rad(M)$ and $\Soc(M)$ are simple functors. By Proposition~\ref{fg}, there is a finite set~$E$ such that $M$ is a quotient, hence a direct summand, of $F=\bigoplusl_{i\in I} k\CC(-,E)$ for some finite set~$I$. Since $k\CC(-,E)\dual\cong k\CC(-,E)$, the dual $M\dual$ is a direct summand of $F\dual\cong F$, and both $M$ and $M\dual$ are generated by their evaluations at $E$.
Thus $M\cong L_{E,M(E)}$ and $M\dual\cong L_{E,M\dual(E)}$, by Lemma~\ref{counit-for-proj}.
As $M$ is a direct summand of $F$ and $M$ is indecomposable, $M$ is a direct summand of~$k\CC(-,E)$, by the Krull-Remak-Schmidt Theorem (Proposition~\ref{KRS}).
So there is a primitive idempotent~$e$ of $k\CC(E,E)\cong\End_{\CF_k}\!\big(k\CC(-,E)\big)$ such that $M\cong k\CC(-,E)e$, and we can assume that $M=k\CC(-,E)e$.\par

If $V$ is a finite dimensional $k$-vector space, and $W$ is a subspace of $V$, set
$$W^\perp=\{\varphi\in \Hom_k(V,k)\mid \varphi(W)=0\}\mpoint$$
If $N$ is a subfunctor of $M$, the assignment sending a finite set $X$ to $N(X)^\perp$ defines a subfunctor $N^\perp$ of $M\dual$, and moreover $N\mapsto N^\perp$ is an order reversing bijection between the set of subfunctors of $M$ and the set of subfunctors of $M\dual$. In particular $\Soc(M)^\perp=\Rad(M\dual)$. Hence $\Soc(M)^\perp(E)=\big(\Soc(M)(E)\big)^\perp=\Rad(M\dual)(E)$.\par

Now $M\dual\neq \Rad(M\dual)$, and $M\dual$ is generated by $M\dual(E)$. Hence $\Rad(M\dual)(E)\neq M\dual(E)$. It follows that $\Soc(M)(E)\neq 0$. Then $\Soc(M)(E)\subseteq \CR_Ee$, and $\Soc(M)(E)$ is a left ideal of~$\CR_E$. It follows that $\Soc(M)(E)$ is not contained in the kernel of the map $t$ defined in Theorem~\ref{self-dual}, that is $t\big(\Soc(M)(E)\big)\neq 0$. Hence
$$0\neq t\big(\Soc(M)(E)\big)=t\big(\Soc(M)(E)e\big)=t\big(e\Soc(M)(E)\big)\mvirg$$
and in particular $e\Soc(M)(E)\neq 0$. Since 
$$e\Soc(M)(E)\cong\Hom_{\CF_k}\big(k\CC(-,E)e,\Soc(M)\big)\mvirg$$
there is a nonzero morphism from $M=k\CC(-,E)e$ to $\Soc(M)$, hence a nonzero morphism from $M/\Rad(M)$ to $\Soc(M)$. Since $M/\Rad(M)$ and $\Soc(M)$ are simple, it is an isomorphism.\mpn

(b) First, by Proposition~\ref{KRS}, both $\Hom_{\CF_k}(M,N)$ and $\Hom_{\CF_k}(N,M)$  are finite dimensional $k$-vectors spaces. \par

Now we can again assume that $M$ is an indecomposable projective and injective functor. For a finitely generated functor $N$, set $\alpha(N)=\dim_k\Hom_{\CF_k}(M,N)$ and $\beta(N)=\dim_k\Hom_{\CF_k}(N,M)$. If $0\to N_1\to N_2\to N_3\to 0$ is a short exact sequence of finitely generated functors, then $\alpha(N_2)=\alpha(N_1)+\alpha(N_3)$ because $M$ is projective, and $\beta(N_2)=\beta(N_1)+\beta(N_3)$ because $M$ is injective. So, in order to prove~(b), as $N$ has finite length, it is enough to assume that $N$ is simple. In that case $\alpha(N)=\dim_k\End_{\CF_k}(N)$ if $M/\Rad(M)\cong N$, and $\alpha(N)=0$ otherwise. Similarly $\beta(N)=\dim_k\End_{\CF_k}(N)$ if $\Soc(M)\cong N$, and $\beta(N)=0$ otherwise. Hence (b) follows from (a).
\endpf


\section{The noetherian case} \label{Section-noetherian}

\bigskip
\noindent
In this section, we shall assume that the ground ring $k$ is noetherian, in which case we obtain more results about subfunctors. For instance, we shall prove that any subfunctor of a finitely generated functor is finitely generated.
It would be interesting to see if the methods developed recently by Sam and Snowden~\cite{SS} for showing noetherian properties of representations of categories can be applied for proving the results of this section.\par

Our first results hold without any assumption on~$k$.

\result{Notation} Let $k$ be a commutative ring, let $E$ be a finite set, and let $M$ be a correspondence functor over $k$. We set
$$\sur{M}(E):=M(E)\Big/\sum_{E'\subset E}k\CC(E,E')M(E')\mvirg$$
where the sum runs over proper subsets $E'$ of $E$.
\fresult

Note that if $F$ is any set of cardinality smaller than $|E|$, then there exists a bijection $\sigma:E'\to F$, where $E'$ is a proper subset of $E$. It follows that $k\CC(E,F)M(F)=k\CC(E,F)R_{\sigma}M(E')\subseteq\sum_{E'\subset E}\limits k\CC(E,E')M(E')$, where $R_\sigma\subseteq F\times E'$ is the graph of~$\sigma$.\par

Note also that $\sur{M}(E)$ is a left module for the essential algebra $\CE_E$, because the ideal
$I_E=\sum\limits_{|Y|<|E|}k\CC(E,Y)k\CC(Y,E)$ of the algebra $\CR_E=k\CC(E,E)$ acts by zero on~$\sur{M}(E)$.

\result{Lemma}\label{localization} Let $k$ be a commutative ring, and let $E$ be a finite set. Let  $M$ be a correspondence functor over~$k$. If $\MF{p}$ is a prime ideal of~$k$, denote by $M_\MF{p}$ the localization of $M$ at~$\MF{p}$, defined by $M_{\MF{p}}(E)=M(E)_{\MF{p}}$ for every finite set~$E$.
\begin{enumerate}
\item $M_\mathfrak{p}$ is a correspondence functor over the localization~$k_\MF{p}$.
\item If $M$ is finitely generated over $k$, then $M_\mathfrak{p}$ is finitely generated over $k_{\MF{p}}$.
\item For each finite set $E$, there is an isomorphism of $k_\mathfrak{p}\CC(E,E)$-modules 
$$\sur{M}(E)_\mathfrak{p}\cong \sur{M_\mathfrak{p}}(E)\mpoint$$
\end{enumerate}
\fresult

\pf (a) This is straightforward.\mpn

(b) If $E$ is a finite set, then clearly $k\CC(-,E)_\mathfrak{p} \cong k_\mathfrak{p}\CC(-,E)$, because this the localization of a free module (on every evaluation).
If $M$ is finitely generated, then there is a finite set $F$ such that $M$ is a quotient of $\bigoplusl_{i\in I}k\CC(-,F)$ for some finite set~$I$.
Then $M_{\MF{p}}$ is a quotient of the functor
$\bigoplusl_{i\in I}k\CC(-,F)_{\MF{p}}\cong \bigoplusl_{i\in I}k_{\MF{p}}\CC(-,F)$,
hence it is a finitely generated functor over~$k_{\MF{p}}$.\mpn

(c) Since localization is an exact functor, the exact sequence of $k$-modules
$$\dirsum{E'\subset E}k\CC(E,E')M(E')\to M(E)\to \sur{M}(E)\to 0$$
gives the exact sequence of $k_{\MF{p}}$-modules 
$$\dirsum{E'\subset E}\big(k\CC(E,E')M(E')\big)_{\MF{p}}\to M(E)_{\MF{p}}\to \sur{M}(E)_{\MF{p}}\to 0\mpoint$$
Now clearly $\big(k\CC(E,E')M(E')\big)_{\MF{p}}=k_{\MF{p}}\CC(E,E')M(E')_{\MF{p}}=k_{\MF{p}}\CC(E,E')M_{\MF{p}}(E')$ for each $E'\subset E$. Hence we get an exact sequence
$$\dirsum{E'\subset E}k_{\MF{p}}\CC(E,E')M_{\MF{p}}(E')\to M_{\MF{p}}(E)\to\sur{M}(E)_{\MF{p}}\to 0\mvirg$$
and it follows that $\sur{M}(E)_{\MF{p}}\cong \sur{M_{\MF{p}}}(E)$.
\endpf

\result{Proposition} \label{localization2} Let $k$ be a commutative ring, let $E$ be a finite set, and let $M$ be a correspondence functor such that $\sur{M}(E)\neq 0$. 
\begin{enumerate}
\item There exists a prime ideal $\mathfrak{p}$ of $k$ such that $\sur{M_\mathfrak{p}}(E)\neq 0$.
\item If moreover $M(E)$ is a finitely generated $k$-module, then there exist subfunctors $A$ and $B$ of $M_{\MF{p}}$ such that $\MF{p}M_{\MF{p}}\subseteq A\subset B$, and a simple module $V$ for the essential algebra $\CE_E$ of $E$ over $k(\MF{p})$ such that $B/A\cong S_{E,V}$, where $k(\MF{p})=k_{\MF{p}}/\MF{p}k_{\MF{p}}$.
\item In this case, there exist positive numbers $c$ and~$d$ such that 
$$c\,|E|^{|X|}\leq \dim_{k(\MF{p})}(M_\MF{p}/\MF{p}M_\MF{p})(X)$$ 
whenever $X$ is a finite set such that $|X|\geq d$.
\end{enumerate}
\fresult

\pf (a) This follows from the well-known fact that the localization map $\sur{M}(E)\longrightarrow\prod_{\MF{p}\in{\rm Spec}(k)}\limits\sur{M}(E)_{\MF{p}}$ is injective, and from the isomorphism $\sur{M}(E)_{\MF{p}}\cong \sur{M_{\MF{p}}}(E)$ of Lemma~\ref{localization}.\mpn

(b) Set $N=M_\MF{p}/\MF{p}M_\MF{p}$ where $\mathfrak{p}$ is the prime ideal obtained in~(a). Then $N$ is a correspondence functor over~$k(\MF{p})$. Suppose that $\sur{N}(E)=0$. Then 
$$M_{\MF{p}}(E)=\MF{p}M_{\MF{p}}(E)+\sum_{E'\subset E}k_{\MF{p}}\CC(E,E')M_{\MF{p}}(E')\mpoint$$
Since $M(E)$ is a finitely generated $k$-module, $M_{\MF{p}}(E)$ is a finitely generated $k_{\MF{p}}$-module, and Nakayama's lemma implies that 
$$M_{\MF{p}}(E)=\sum_{E'\subset E}\limits k_{\MF{p}}\CC(E,E')M_{\MF{p}}(E')\mvirg$$
that is, $\sur{M_{\MF{p}}}(E)=0$. This contradicts (a) and shows that $\sur{N}(E)\neq 0$.\par

Now $\sur{N}(E)$ is a nonzero module for the essential algebra $\CE_E$ of $E$ over $k(\MF{p})$, and it is finite dimensional over $k(\MF{p})$ (because $M_{\MF{p}}(E)$ is a finitely generated $k_{\MF{p}}$-module). Hence it admits a simple quotient $V$ as $\CE_E$-module. Then $V$ can be viewed as a simple $k(\MF{p})\CR_E$-module by inflation, and it is also a quotient of~$N(E)$.
By Proposition~\ref{subquotients}, there exist subfunctors $A/\MF{p}M_{\MF{p}}\subset B/\MF{p}M_{\MF{p}}$ of~$N$ such that $B/A$ is isomorphic to the simple functor $S_{E,V}$, proving~(b). \mpn

(c) By (b) and Theorem~\ref{exponential-simple}, there exist positive numbers $c$ and $d$ such that 
$$c\,|E|^{|X|}\leq \dim_{k(\MF{p})}(B/A)(X) $$ 
whenever $X$ is a finite set such that $|X|\geq d$. Assertion (c) follows.
\endpf

Now we assume that $k$ is noetherian and we can state the critical result.

\result{Theorem} \label{bounded size} Let $k$ be a commutative noetherian ring.
Let $N$ be a subfunctor of a correspondence functor $M$ over~$k$.
If $E$ and $Y$ are finite sets such that $M$ is generated by $M(E)$ and $\sur{N}(Y)\neq 0$, then $|Y|\leq 2^{|E|}$.
\fresult

\pf  Since $M$ is generated by $M(E)$, choosing a set $I$ of generators of $M(E)$ yields a surjection
$\Phi: P=\dirsum{i\in I}k\CC(-,E)\to M$. Let $L=\Phi^{-1}(N)$. Since $\Phi$ induces a surjection $\sur{L}(Y)\to \sur{N}(Y)$, and since $P$ is generated by $P(E)$, we can replace $M$ by $P$ and $N$ by~$L$.
Hence we now assume that $N$ is a subfunctor of $\dirsum{i\in I}k\CC(-,E)$.\par

Since $\sur{N}(Y)\neq 0$, there exists
$$m\in N(Y)-\sum_{Y'\subset Y}\limits k\CC(Y,Y')N(Y') \mpoint$$
Let $N'$ be the subfunctor of $N$ generated by~$m$. Then clearly $\sur{N'}(Y)\neq 0$, because
$$m\in N'(Y)-\sum_{Y'\subset Y}\limits k\CC(Y,Y')N'(Y') \mpoint$$
Moreover $N'(Y)=k\CC(Y,Y)m$ is a finitely generated $k$-module, and there is a finite subset $S$ of $I$ such that $m\in \dirsum{i\in S}k\CC(Y,E)$. Therefore $N'\subseteq \dirsum{i\in S}k\CC(-,E)$.
Replacing $N$ by $N'$, we can assume moreover that the set $I$ is finite. In other words, there exists an integer $s\in\N$ such that $N\subseteq k\CC(-,E)^{\oplus s}$.\par

Now by Proposition~\ref{localization2}, there exists a prime ideal $\MF{p}$ of $k$ such that $\sur{N_{\MF{p}}}(Y)\neq 0$. Moreover $N(Y)$ is a submodule of $k\CC(Y,E)^{\oplus s}$, which is a finitely generated (free) $k$-module. Since $k$ is noetherian, it follows that $N(Y)$ is a finitely generated $k$-module.\par

By Proposition~\ref{localization2}, there exist subfunctors $A\subset B$ of $N_{\MF{p}}$ such that $B/A$ is isomorphic to a simple functor of the form $S_{Y,V}$, where $V$ is a simple module for the essential algebra of $Y$ over~$k(\MF{p})$. In particular $Y$ is minimal such that $(B/A)(Y)\neq 0$, thus $(\sur{B/A})(Y)\cong (B/A)(Y)\cong V$.\par

It follows that $\sur{B}(Y)\neq 0$, and $B$ is a subfunctor of $k_{\MF{p}}\CC(-,E)^{\oplus s}$.
In other words, replacing $k$ by $k_{\MF{p}}$ and $N$ by $B$, we can assume that $k$ is a noetherian local ring, that $\MF{p}$ is the unique maximal ideal of~$k$, and that $N$ has a subfunctor $A$ such that $N/A$ is isomorphic to $S_{Y,V}$, where $V$ is a simple module for the essential algebra~$\CE_Y$ over~$k/\MF{p}$.\par

We claim that there exists an integer $n\in\N$ such that
$$N(Y)\neq A(Y)+\big(\MF{p}^n\CC(Y,E)^{\oplus s}\cap N(Y)\big)\mpoint$$
Indeed $N(Y)$ is a submodule of the finitely generated $k$-module $k\CC(Y,E)^{\oplus s}$.
By the Artin-Rees lemma (see Theorem 8.5 in~\cite{Ma}), there exists an integer $l\in \N$ such that for any $n\geq l$
$$\MF{p}^n\CC(Y,E)^{\oplus s}\cap N(Y)=\MF{p}^{n-l}\big(\MF{p}^l\CC(Y,E)^{\oplus s}\cap N(Y)\big)\mpoint$$
Let $m_1,\ldots ,m_r$ be generators of $N(Y)$ as a $k$-module. Suppose that $n>l$ and that
$N(Y)= A(Y)+\big(\MF{p}^n\CC(Y,E)^{\oplus s}\cap N(Y)\big)$. Then
$$N(Y)= A(Y)+\MF{p}^{n-l}\big(\MF{p}^l\CC(Y,E)^{\oplus s}\cap N(Y)\big) \mpoint$$
It follows that for each $i\in\{1,\ldots,r\}$, there exist $a_i\in A(Y)$ and scalars $\lambda_{i,j}\in \MF{p}^{n-l}$, for $1\leq j\leq r$, such that
$$m_i=a_i+\sum_{j=1}^r\lambda_{i,j}m_j\mpoint$$
In other words the sequence $(a_i)_{i=1,\ldots,r}$ is the image of the sequence $(m_i)_{i=1,\ldots,r}$ under the matrix $J=\Id_r-\Lambda$, where $\Lambda$ is the matrix of coefficients $\lambda_{i,j}$, and $\Id_r$ is the identity matrix of size $r$. Since $\Lambda$ has coefficients in $\MF{p}^{n-l}\subseteq \MF{p}$, 
the determinant of $J$ is congruent to 1 modulo $\MF{p}$, hence $J$ is invertible. It follows that the $m_i$'s are linear combinations of the $a_j$'s with coefficients in $k$. Hence $m_i\in A(Y)$ for $1\leq i\leq r$, thus $N(Y)=A(Y)$.
This is a contradiction since $(N/A)(Y)\cong V\neq 0$. This proves our claim.\par

We have obtained that $N\neq A+\big(\MF{p}^n\CC(-,E)^{\oplus s}\cap N\big)$. Since $N/A$ is simple, it follows that $\MF{p}^n\CC(-,E)^{\oplus s}\cap N=\MF{p}^n\CC(-,E)^{\oplus s}\cap A$. \par

Now we reduce modulo $\MF{p}^n$ and we let respectively $\hat{A}$ and $\hat{N}$ denote the images of $A$ and $N$ in the reduction $(k/\MF{p}^n)\CC(-,E)^{\oplus s}$. Then 
\begin{eqnarray*}
\hat{N}/\hat{A}\!\!\!&=&\!\!\!\Big(\big(N+\MF{p}^n\CC(-,E)^{\oplus s}\big)/\MF{p}^n\CC(-,E)^{\oplus s}\Big)
\Big/\Big(\big(A+\MF{p}^n\CC(-,E)^{\oplus s}\big)/\MF{p}^n\CC(-,E)^{\oplus s}\Big)\\
&\cong&\Big(N/\big(\MF{p}^n\CC(-,E)^{\oplus s}\cap N\big)\Big)\Big/\Big(A/\big(\MF{p}^n\CC(-,E)^{\oplus s}\cap A\big)\Big)\\
&\cong&N/A\mvirg
\end{eqnarray*}
and this is isomorphic to the simple functor $S_{Y,V}$ over the field $k/\MF{p}$.
Hence for any finite set $X$, the module $\hat{N}(X)/\hat{A}(X)$ is a $(k/\MF{p})$-vector space. Moreover, by Proposition~\ref{localization2}, there exist positive numbers $c$ and $d$ such that the dimension of this vector space is larger than $c\,|Y|^{|X|}$ whenever $|X|\geq d$.\par

Now for any finite set $X$, the module $(k/\MF{p}^n)\CC(X,E)^{\oplus s}$ is filtered by the submodules $\Gamma_j=(\MF{p}^j/\MF{p}^n)\CC(X,E)^{\oplus s}$, for $j=1,\ldots,n-1$, and the quotient $\Gamma_j/\Gamma_{j+1}$ is a vector space over $k/\MF{p}$, of dimension $sd_j\,2^{|X||E|}$, where $d_j=\dim_{k/\MF{p}}(\MF{p}^j/\MF{p}^{j+1})$.
It follows that, for $|X|\geq d$,
$$c\,|Y|^{|X|}\leq\dim_{k/\MF{p}}\big(\hat{N}(X)/\hat{A}(X)\big)\leq s\big(\sum_{j=0}^{n-1}d_j\big)2^{|X||E|}\mpoint$$
As $|X|$ tends to infinity, this forces $|Y|\leq 2^{|E|}$, completing the proof of Theorem~\ref{bounded size}.
\endpf

\result{Corollary} \label{N subseteq M} Let $k$ be a commutative noetherian ring and let $N$ be a subfunctor of a correspondence functor $M$ over~$k$.
\begin{enumerate}
\item If $E$ is a finite set such that $M$ is generated by $M(E)$ and if $F$ is a finite set with $|F|\geq 2^{|E|}$, then $N$ is generated by $N(F)$.
\item If $M$ has bounded type, then $N$ has bounded type. In particular, over~$k$, any correspondence functor of bounded type has a bounded presentation. 
\item If $M$ is finitely generated, then $N$ is finitely generated. In particular, over~$k$, any finitely generated correspondence functor is finitely presented.
\end{enumerate}
\fresult

\pf
(a) Let $E$ be a finite set such that $M$ is generated by~$M(E)$. If $X$ is a finite set such that $\sur{N}(X)\neq 0$, then $|X|\leq 2^{|E|}$, by Theorem~\ref{bounded size}. For each integer $e\leq 2^{|E|}$, let $[e]=\{1,\ldots,e\}$ and choose a subset $S_e$ of $N([e])$ which maps to a generating set of $\sur{N}([e])$ as a $k$-module.
Each $i\in S_e$ yields a morphism $\psi_{e,i}:k\CC(-,[e])\to N$. Let
$$Q=\dirsum{e\leq 2^{|E|}}\dirsum{i\in S_e}k\CC(-,[e]) \qquad\text{and}\qquad
\Psi=\sumb{e\leq 2^{|E|}}{i\in S_e}\psi_{e,i}: Q\to N\mpoint$$
Then by construction the induced map
$$\sur{\Psi}_X:\sur{Q}(X)\to\sur{N}(X)$$
is surjective, for any finite set $X$, because either $\sur{N}(X)=0$ or $|X|=e\leq 2^{|E|}$.
Suppose that $\Psi:Q\to N$ is not surjective and
let $A$ be a set of minimal cardinality such that $\Psi_A:Q(A)\to N(A)$ is not surjective. Let $l\in N(A)-\Im\Psi_A$.
Since the map $\sur{\Psi}_A$ is surjective, there is an element $q\in Q(A)$ and elements $l_e\in N([e])$ and $R_e\in k\CC(A,[e])$, for $e<|A|$, such that 
$$l=\Psi_A(q)+\sum_{e<|A|}R_el_e\mpoint$$
The minimality of $A$ implies that the map $\Psi_{[e]}:Q([e])\to N([e])$ is surjective for each $e<|A|$, so there are elements $q_e\in Q([e])$, for $e<|A|$, such that $\Psi_{[e]}(q_e)=l_e$. It follows that $l=\Psi_A\big(q+\sum\limits_{e<|A|}R_eq_e\big)$, thus $l\in\Im\Psi_A$. This contradiction proves that the morphism $\Psi:Q\to N$ is surjective.\par

Now let $F$ be a set with $|F|\geq 2^{|E|}$. For each $e\leq 2^{|E|}$, the representable functor $k\CC(-,[e])$ is generated by its evaluation at~$[e]$, hence also by its evaluation at~$F$, because $k\CC(-,[e])$ is a direct summand of~$k\CC(-,F)$ by Corollary~\ref{summand2}. Therefore $Q$ is generated by $Q(F)$.
Since $\Psi:Q\to N$ is surjective, it follows that $N$ is generated by $N(F)$.\mpn

(b) This follows clearly from (a).\mpn

(c) If now $M$ is finitely generated, then the same argument applies,
but we can assume moreover that all the sets $S_e$ appearing in the proof of (a) are {\em finite}, since for any finite set $X$,
the module $N(X)$ is finitely generated, being a submodule of the finitely generated module $M(X)$.
It follows that the functor $Q$ of the proof of~(a) is finitely generated and this proves~(c).
\endpf

It follows from (b) and Theorem~\ref{counit-iso} that, whenever $k$ is noetherian, any correspondence functor of bounded type is isomorphic to $L_{F,V}$ for some $F$ and~$V$. We shall return to this in Theorem~\ref{All-SEV} below.

\result{Notation} We denote by $\CF_k^b$ the full subcategory of $\CF_k$ consisting of correspondence functors having bounded type and by $\CF_k^f$ the full subcategory of $\CF_k$ consisting of finitely generated functors.
\fresult

\result{Corollary} \label{quasi-finitely abelian} Let $k$ be a commutative noetherian ring. Then the categories $\CF_k^b$ and $\CF_k^f$ are abelian full subcategories of $\CF_k$.
\fresult

\pf Any quotient of a functor of bounded type has bounded type and any quotient of a finitely generated functor is finitely generated. When $k$ is noetherian, any subfunctor of a functor of bounded type has bounded type and any subfunctor of a finitely generated functor is finitely generated, by Corollary~\ref{N subseteq M}.
\endpf


\bigbreak
\section{Stabilization results} \label{Section-stabilization}

\bigskip
\noindent
Recall from Lemma~\ref{JXW} that for any finite set $E$ and any $\CR_E$-module $V$, we have defined a subfunctor $J_{E,V}$ of $L_{E,V}$ by setting
$$J_{E,V}(X)=\Big\{\sum_i R_i\mathop{\otimes}_{\CR_E}\limits v_i \mid R_i\in \CC(X,E),\;v_i\in V,\; \forall S\in\CC(E,X),\;\sum_i(SR_i)v_i=0 \Big\}\mpoint$$
Moreover $J_{E,V}(E)=0$ and $S_{E,V}=L_{E,V}/J_{E,V}$.\par

We have seen in Proposition~\ref{induction} that $L_{F,V\!\uparrow_E^F}\cong L_{E,V}$ whenever $|F|\geq |E|$.
The subfunctor $J_{F,V\!\uparrow_E^F}$ vanishes at~$F$, hence also at~$E$, so that $J_{F,V\!\uparrow_E^F}\subset J_{E,V}$.
When $k$ is noetherian, we show that this decreasing sequence reaches zero.

\result{Theorem} \label{JEV-zero} Let $k$ be a commutative noetherian ring, let $E$ be a finite set, and let $V$ be an $\CR_E$-module. For any finite set $F$ such that $|F|\geq 2^{|E|}$, we have $J_{F,V\!\uparrow_E^F}=0$.
\fresult

\pf Let $F$ be a finite set. By Proposition~\ref{induction}, there is an isomorphism $L_{F,V\!\uparrow_E^F}\cong L_{E,V}$. Thus $J_{F,V\!\uparrow_E^F}$ is isomorphic to a subfunctor of~$L_{E,V}$. Since $J_{F,V\!\uparrow_E^F}(F)=0$, it follows from Corollary~\ref{vanishing-below} that $J_{F,V\!\uparrow_E^F}(X)=0$ for any finite set $X$ with $|X|\leq |F|$.\par

We now assume that $J_{F,V\!\uparrow_E^F}\neq 0$ and we prove that $|F|< 2^{|E|}$.
Let $Y$ be a set of minimal cardinality such that $J_{F,V\!\uparrow_E^F}(Y)\neq 0$. Then $|Y|>|F|$.
Moreover $\sur{J_{F,V\!\uparrow_E^F}}(Y)\cong J_{F,V\!\uparrow_E^F}(Y) \neq 0$, hence $|Y|\leq 2^{|E|}$ by Theorem~\ref{bounded size}, because $J_{F,V\!\uparrow_E^F}$ is (isomorphic to) a subfunctor of $L_{E,V}$, which is generated by $L_{E,V}(E)=V$. It follows that $|F|<|Y|\leq 2^{|E|}$.
\endpf

We now show that, over a noetherian ring, any correspondence functor of bounded type is isomorphic to $L_{F,V}$ for some $F$ and~$V$, or also isomorphic to $S_{G,W}$ for some $G$ and~$W$ (where the symbol $S_{G,W}$ refers to Notation~\ref{SXW}).

\result{Theorem} \label{All-SEV} Let $k$ be a commutative noetherian ring.
Let $M$ be a correspondence functor over $k$ generated by $M(E)$, for some finite set~$E$. 
\begin{enumerate}
\item For any finite set $F$ such that $|F|\geq 2^{{|E|}}$, the counit morphism $\eta_{M,F}:L_{F,M(F)}\to M$ is an isomorphism.
\item For any finite set $G$ such that $|G|\geq 2^{2^{{|E|}}}$, we have $M\cong L_{G,M(G)}$ and $J_{G,M(G)}=0$, hence $M\cong S_{G,M(G)}$.
\end{enumerate}
\fresult

\pf
(a) If $M$ is generated by $M(E)$, then there is a set $I$ and a surjective morphism $P=\dirsum{i\in I}k\CC(-,E)\to M$. If $F$ is a finite set with $|F|\geq 2^{|E|}$, then by Corollary~\ref{N subseteq M} the kernel $K$ of this morphism is generated by $K(F)$. Then $K$ is in turn covered by a projective functor~$Q$ and we have a bounded presentation
$$Q\to P\to M\to 0$$
with both $Q$ and $P$ generated by their evaluation at~$F$. By Theorem~\ref{counit-iso}, the counit morphism $\eta_{M,F}:L_{F,M(F)}\to M$ is an isomorphism.\mpn

(b) For any finite set $F$ such that $|F|\geq 2^{{|E|}}$, we have $M\cong L_{F,M(F)}$ by~(a).
For any finite set $G$ such that $|G|\geq 2^{|F|}$, that is, $|G|\geq 2^{2^{{|E|}}}$, we obtain $J_{G,M(F)\!\uparrow_F^G}=0$ by Theorem~\ref{JEV-zero}.
It follows that
$$M\cong L_{F,M(F)} \cong L_{G,M(F)\!\uparrow_F^G} = S_{G,M(F)\!\uparrow_F^G} \mpoint$$
Finally, notice that, by the definition of~$L_{F,M(F)}$, we have $M(G)\cong L_{F,M(F)}(G)=M(F)\!\uparrow_F^G$,
so we obtain $M\cong L_{G,M(G)} = S_{G,M(G)}$.
\endpf

Other kinds of stabilizations also occur, as the next theorems show.

\result{Theorem} \label{stabilization} Let $k$ be a commutative noetherian ring, let $M$ and $N$ be correspondence functors over $k$, and let $E$ and $F$ be finite sets.
\begin{enumerate}
\item If $M$ is generated by $M(E)$, then for $|F|\geq 2^{|E|}$, the evaluation map at $F$
$$\Hom_{\CF_k}(M,N)\to\Hom_{\CR_F}\!\big(M(F),N(F)\big)$$
is an isomorphism.
\item If $M$ has bounded type, then for any integer $i\in \N$, there exists an integer~$n_i$ such that the evaluation map
$$\Ext_{\CF_k}^i(M,N)\to \Ext_{\CR_F}^i\!\big(M(F),N(F)\big)$$
is an isomorphism whenever $|F|\geq n_i$.
\end{enumerate}
\fresult

\pf (a) By Theorem~\ref{All-SEV}, we have an isomorphism $M\cong L_{F,M(F)}$ for $|F|\geq 2^{|E|}$. Hence
$$\Hom_{\CF_k}(M,N)\cong \Hom_{\CF_k}\!\big(L_{F,M(F)},N\big)\cong\Hom_{\CR_F}\!\big(M(F),N(F)\big)\mvirg$$
where the last isomorphism comes from the adjunction of Lemma~\ref{adjunction}, and is given by evaluation at~$F$.\mpn

(b) This assertion will follow from (a) by {\em d\'ecalage} and induction on~$i$.
If $M$ is generated by $M(E)$, then there is an exact sequence of correspondence functors
$$0\to L\to P\to M\to 0$$
where $P$ is projective and generated by $P(E)$. This gives an exact sequence
$$0\to \Hom_{\CF_k}(M,N)\to \Hom_{\CF_k}(P,N)\to \Hom_{\CF_k}(L,N)\to\Ext_{\CF_k}^1(M,N)\to 0\mvirg$$
and isomorphisms $\Ext_{\CF_k}^i(M,N)\cong \Ext_{\CF_k}^{i-1}(L,N)$ for $i\geq 2$.\par
Now $L$  has bounded type by Corollary~\ref{N subseteq M}, and $P(F)$ is a projective $\CR_F$-module by Lemma~\ref{counit-for-proj}, whenever $|F|$ is large enough. It follows that there is also an exact sequence
\begin{eqnarray*}
0\to\!\Hom_{\CR_F}\!\big(M(F),N(F)\big)\!\!\!\!&\to&\!\!\!\! \Hom_{\CR_F}\!\big(P(F),N(F)\big)\to\\
&&\!\!\!\! \Hom_{\CR_F}\!\big(L(F),N(F)\big)\to\!\Ext_{\CR_F}^1\!\big(M(F),N(F)\big)\to\! 0
\end{eqnarray*}
and isomorphisms $\Ext_{\CR_F}^i\!\big(M(F),N(F)\big)\cong\Ext_{\CR_F}^{i-1}\!\big(L(F),N(F)\big)$ for $i\geq 2$, whenever $F$ is large enough.\par

Now by (a), the exact sequences
$$0\to \Hom_{\CF_k}(M,N)\to \Hom_{\CF_k}(P,N)\to \Hom_{\CF_k}(L,N)$$
and 
$$0\to\Hom_{\CR_F}\!\big(M(F),N(F)\big)\to\Hom_{\CR_F}\!\big(P(F),N(F)\big)\to
\Hom_{\CR_F}\!\big(L(F),N(F)\big)$$
are isomorphic for $F$ large enough. It follows that
$$\Ext_{\CF_k}^1(M,N)\cong \Ext_{\CR_F}^1\!\big(M(F),N(F)\big)
\mpoint$$
Similarly, for each $i\geq 2$, when $F$ is large enough (depending on $i$), there are isomorphisms $\Ext_{\CF_k}^i(M,N)\cong\Ext_{\CR_F}^i\!\big(M(F),N(F)\big)$.
\endpf

There is also a stabilization result involving the $\Tor$ groups.
 
\result{Theorem} Let $k$ be a commutative noetherian ring, and $E$ be a finite set. If $F$ is a finite set with $|F|\geq 2^{2^{|E|}}$, then for any finite set $X$ and any left $\CR_E$-module~$V$, we have
$$\Tor_1^{\CR_F}\!\big(k\CC(X,F),V\!\uparrow_E^F\big)=0\mpoint$$
\vspace{-2ex}
\fresult

\pf Let $V$ be a left $\CR_E$-module and $s:Q\to V$ be a surjective morphism of $\CR_E$-modules, where $Q$ is projective. Let $K$ denote the kernel of the surjective morphism
$$L_{E,s}:L_{E,Q}\to L_{E,V}\mpoint$$
Since $L_{E,Q}$ is generated by $L_{E,Q}(E)\cong Q$, it follows from Corollary~\ref{N subseteq M} that $K$ is generated by $K(G)$ whenever $G$ is a finite set with $|G|\geq 2^{|E|}$.
Now by Theorem~\ref{All-SEV}, the counit $L_{F,K(F)}\to K$ is an isomorphism whenever $F$ is a finite set with $|F|\geq 2^{|G|}$. Hence if $|F|\geq 2^{2^{|E|}}$,
we have an exact sequence of correspondence functors
\begin{equation}\label{exact sequence of L}
0\to L_{F,W_F}\to L_{E,Q}\to L_{E,V}\to 0\mvirg
\end{equation}
where $W_F=K(F)$, and where the middle term $L_{E,Q}$ is projective. Evaluating this sequence at $F$, we get the exact sequence of $\CR_F$-modules
$$0\to W_F\to k\CC(F,E)\otimes_{\CR_E}Q\to k\CC(F,E)\otimes_{\CR_E}V\to 0 \mvirg$$
where the middle term is projective.\par

Let $X$ be a finite set. Applying the functor $k\CC(X,F)\otimes_{\CR_F}(-)$ to this sequence yields the exact sequence
\begin{multline*}0\to\Tor_1^{\CR_F}\!\big(k\CC(X,F),V\!\uparrow_E^F\big)\to \\k\CC(X,F)\otimes_{\CR_F}W_F\to k\CC(X,E)\otimes_{\CR_E}Q\to k\CC(X,E)\otimes_{\CR_E}V\to 0\mvirg
\end{multline*}
because $k\CC(X,F)\otimes_{\CR_F}k\CC(F,E)\cong k\CC(X,E)$ by Corollary~\ref{composition-iso}, as $|F|\geq |E|$.
On the other hand, evaluating at $X$ the exact sequence~\ref{exact sequence of L} gives the exact sequence
$$0\to  k\CC(X,F)\otimes_{\CR_F}W_F\to k\CC(X,E)\otimes_{\CR_E}Q\to k\CC(X,E)\otimes_{\CR_E}V\to 0\mpoint$$
In both latter exact sequences, the maps are exactly the same. It follows that
$$\Tor_1^{\CR_F}\!\big(k\CC(X,F),V\!\uparrow_E^F\big)=0 \mvirg$$
as was to be shown.
\endpf

As a final approach to stabilization, we introduce the following definition.

\result{Definition} Let $\CG_k$ denote the following category:
\begin{itemize}
\item The objects of $\CG_k$ are pairs $(E,U)$ consisting of a finite set $E$ and a left $\CR_E$-module $U$.
\item A morphism $\varphi:(E,U) \to (F,V)$ in $\CG_k$ is a morphism of $\CR_E$-modules $U\to k\CC(E,F)\otimes_{\CR_F}V$.
\item The composition of morphisms $\varphi:(E,U)\to (F,V)$ and $\psi:(F,V)\to (G,W)$ is the morphism obtained by composition
$$U\stackrel{\varphi}{\to}k\CC(E,F)\otimes_{\CR_F}V\stackrel{\Id\otimes \psi}{\longrightarrow}k\CC(E,F)\otimes_{\CR_F}k\CC(F,G)\otimes_{\CR_G}W\stackrel{\mu\otimes\Id}{\longrightarrow}k\CC(E,G)\otimes_{\CR_G}W\mvirg$$
where $\mu:k\CC(E,F)\otimes_{\CR_F}k\CC(F,G)\to k\CC(E,G)$ is the composition in the category $k\CC$.
\item The identity morphism of $(E,U)$ is the canonical isomorphism 
$$U\to k\CC(E,E)\otimes_{\CR_E}U$$
resulting from the definition $\CR_E=k\CC(E,E)$.
\end{itemize}
\fresult

One can check easily that $\CG_k$ is a $k$-linear category.

\result{Theorem} Let $k$ be a commutative ring. Let $\LL:\CG_k\to\CF_k^b$ be the assignment sending $(E,U)$ to $L_{E,U}$, and $\varphi:(E,U)\to (F,V)$ to the morphism $L_{E,U}\to L_{F,V}$ associated by adjunction to $\varphi:U\to L_{F,V}(E)$.
\begin{enumerate}
\item $\LL$ is a fully faithful $k$-linear functor.
\item $\LL$ is an equivalence of categories if $k$ is noetherian.
\item $\CG_k$ is an abelian category if $k$ is noetherian.
\end{enumerate}
\fresult

\pf It is straightforward to check that $\LL$ is a $k$-linear functor. It is moreover fully faithful, since
$$\Hom_{\CF_k^b}(L_{E,U},L_{F,V})\cong \Hom_{\CR_E}\!\big(U,k\CC(E,F)\otimes_FV\big)\cong\Hom_{\CG_k}\!\big((E,U),(F,V)\big)\mpoint$$
Finally, if $k$ is noetherian, then any correspondence functor $M$ of bounded type is isomorphic to a functor of the form $L_{E,U}$, by Theorem~\ref{All-SEV}, for $E$ large enough and $U=M(E)$. Hence $\LL$ is essentially surjective, so it is an equivalence of categories. In particular, $\CG_k$ is abelian by Corollary~\ref{quasi-finitely abelian}.\endpf

\bigskip
\noindent
Serge Bouc, CNRS-LAMFA, Universit\'e de Picardie - Jules Verne,\\
33, rue St Leu, F-80039 Amiens Cedex~1, France.\\
{\tt serge.bouc@u-picardie.fr}

\medskip
\noindent
Jacques Th\'evenaz, Section de math\'ematiques, EPFL, \\
Station~8, CH-1015 Lausanne, Switzerland.\\
{\tt Jacques.Thevenaz@epfl.ch}

\end{document}